\documentclass[10pt]{amsart}
\usepackage{amssymb}
\usepackage{amsfonts}
\usepackage{graphicx}
\usepackage{amsmath}
\usepackage{accents}
\usepackage{calc}
\usepackage{tikz}

\usetikzlibrary{decorations.markings}
\tikzstyle{vertex}=[circle, draw, inner sep=0pt, minimum size=6pt]

\theoremstyle{plain}

\newtheorem{problem}{Problem}

\numberwithin{equation}{section}
\input{tcilatex}

\begin{document}
	\leftline{\small Scientiae Mathematicae Japonicae (to appear) }
	\vspace{13mm} 
	
\title[]{Groupoid Factorizations in the Semigroup of Binary Systems}
\author{Hiba F. Fayoumi}
\address[H. F. Fayoumi]{University of Toledo\\
2801 Bancroft Street, Toledo, Ohio, 43606, U.S.A. }
\email{hiba.fayoumi@UToledo.edu}
\urladdr{}
\thanks{}
\date{}
\subjclass[2010]{Primary 20N02}
\keywords{groupoid factorization, groupoid decomposition, composite
groupoid, prime groupoid, Bin(X), idempotent groupoid, signature-factor,
similar-factor, orient-factor, skew-factor, u-normal, j-normal, $\Psi $%
-type-factor, $\tau $-type-factor}
\dedicatory{}

\begin{abstract}
Let $(X,\bullet )$ be a groupoid (binary algebra) and $Bin(X\dot{)}$ denote
the collection of all groupoids defined on $X$. We introduce two methods of
factorization for this binary system under the binary groupoid product
\textquotedblleft $\diamond $\textquotedblright\ in the semigroup $\left(
Bin\left( X\right) ,\diamond \right) $. We conclude that a strong
non-idempotent groupoid can be represented as a product of its \textit{%
similar-} and \textit{signature-} derived factors. Moreover, we show that a
groupoid with the orientation property is a product of its \textit{orient-}
and \textit{skew-} factors. These unique factorizations can be useful for
various applications in other areas of study. Application to algebras such
as $B/BCH/BCI/BCK/BH/BI/d$-algebra are widely given throughout this paper.
\end{abstract}

\maketitle

\section{Introduction}

Algebraic structures play a vital role in mathematical applications such as
information science, network engineering, computer science, cell biology,
etc. This encourages sufficient motivation to study abstract algebraic
concepts and review previously obtained results. One such concept of
interest to many mathematicians over the past two decades or so is that of a
simple yet very interesting notion of a single set with one binary
operation, historically known as magma and more recently referred to as
groupoid. Bruck \cite{Br58} published the book,\textquotedblleft A Survey of
Binary Systems\textquotedblright\ in which the theory of groupoids, loops,
quasigroups, and several algebraic structures were discussed. Bor$\mathrm{{%
\dot{u}}}$vka in \cite{Bor76} explained the foundations for the theory of
groupoids, set decompositions and their application to binary systems.

Given a binary operation \textquotedblleft $\bullet $\textquotedblright\ on
a non-empty set $X$, the groupoid $\left( X,\bullet \right) $ is a
generalization of the very well-known structure of a group. H. S. Kim and J.
Neggers in \cite{KN08} investigated the structure $\left( Bin\left( X\right)
,\diamond \right) $ where $Bin\left( X\right) $ is the collection of all
binary systems (groupoids or algebras) defined on a non-empty set $X$ along
with an associative binary product $(X,\ast )\diamond (X,\circ )=\left(
X,\bullet \right) $ such that $x\bullet y=(x\ast y)\circ (y\ast x)$ for all $%
x$, $y\in X$. They recognized that the left-zero-semigroup serves as the
identity of this semigroup. The present author in \cite{Fa11} introduced the
notion of the center $ZBin(X)$ in the semigroup $\left( Bin\left( X\right)
,\diamond \right) $, and proved that $(X,\bullet )\in ZBin(X)$, if and only
if $(X,\bullet )$ is locally-zero. Han and Kim in \cite{HKN12} introduced
the notion of hypergroupoids $HBin(X)$, and showed that $(HBin(X),\diamond )$
is a supersemigroup of the semigroup $(Bin(X),\diamond )$ via the
identification $x\longleftrightarrow \{x\}$. They proved that $(HBin^{\ast
}(X),\ominus ,[\emptyset ])$ is a $BCK$-algebra.

In this paper, we investigate the following problem:

\begin{description}
\item[Main Problem] Consider the semigroup $\left( Bin\left( X\right)
,\diamond \right) $. Let the left-zero-semigroup be denoted as $%
id_{Bin\left( X\right) }$. Given a groupoid (binary system) $(X,\bullet )\in 
$ $Bin(X)$, is it possible to find two groupoid-factors $\left( X,\ast
\right) $ and $\left( X,\circ \right) $ such that 
\begin{equation*}
\left( X,\bullet \right) =(X,\ast )\diamond (X,\circ )\text{?}
\end{equation*}%
If so,

\begin{problem}[Uniqueness]
Are the corresponding groupoid-factors:
\end{problem}

\begin{enumerate}
\item \label{prob1.1}\textit{Distinct, i.e., }$(X,\ast )\neq (X,\circ )$%
\textit{?}

\item \label{prob1.2}\textit{Unique, i.e., if }$\left( X,\bullet \right)
=(X,\ast )\diamond (X,\circ )$\textit{, is it possible for }$\left(
X,\bullet \right) =\left( X,\lhd \right) \diamond $\textit{\ }$\left( X,\rhd
\right) $\textit{\ such that }$(X,\ast )\neq (X,\lhd )$\textit{\ and }$%
(X,\circ )\neq (X,\rhd )$\textit{?}

\item \label{prob1.3}\textit{Different from }$(X,\bullet )$\textit{, i.e., }$%
(X,\ast )\neq (X,\bullet )$\textit{\ and }$(X,\circ )\neq (X,\bullet )$%
\textit{?}

\item \label{prob1.4}\textit{Different from the left-zero-semigroup, i.e., }$%
(X,\ast )\neq id_{Bin\left( X\right) }$\textit{\ and }$(X,\circ )\neq
id_{Bin\left( X\right) }$\textit{?}
\end{enumerate}

\begin{problem}[Derivation]
How do we find the groupoid-factors? Are they:
\end{problem}

\begin{enumerate}
\item \label{prob2.1}\textit{Derived (related to, based off of, dependent
on) from: the parent groupoid }$(X,\bullet )$\textit{?}

\item \label{prob2.2}\textit{Derived from the identity }$id_{Bin\left(
X\right) }$\textit{?}
\end{enumerate}

\begin{problem}[Factorization]
If we use a certain method to find the two groupoid-factors, what is the
nature of this factorization?
\end{problem}

\begin{enumerate}
\item \label{prob3.1}\textit{Is it unique? }

\item \label{prob3.2}\textit{When is it commutative?}
\end{enumerate}
\end{description}

We begin answering these questions by introducing two methods for factoring
a random groupoid in $Bin(X\dot{)}$ using the product \textquotedblleft $%
\diamond $\textquotedblright . We will show that both methods result in
unique factorizations (Problem 3.\ref{prob3.1}) of a given groupoid and
hence we answer Problem 1.\ref{prob1.2} with a definite yes! Section two
provides some definitions and preliminary ideas which are necessary in this
context. We also present a summarized table of \textquotedblleft
logic\textquotedblright\ algebras for a clear view. Section three describes $%
AU$- and $UA$-factorizations, which comprises the first method (method-1) of
factoring. In fact, method-1 factors a groupoid $\left( X,\bullet \right) $
by obtaining two \textit{derived} factors from it (Problem 2.\ref{prob2.1})
and from the left-zero-semigroup (Problem 2.\ref{prob2.2}), the \textit{%
signature-} and \textit{similar-}factors, respectively. We prove that a
strong groupoid has a commutative method-1 factorization (Problem 3.\ref%
{prob3.2}). The possibility of this first method is shown to be feasible and
produces non-trivial decompositions (Problem 1.\ref{prob1.4}), however, it
is restricted to non-idempotent groupoids only. Hence, section four
introduces an $OJ$\textit{-} and a $JO$\textit{-}factorization, which
constitutes our second method (method-2). We will demonstrate that the
latter method is sufficient for idempotent as well as non-idempotent
groupoids. In addition, an interesting outcome of method-2 is that one of
the factors is not \textit{derived} from the parent groupoid (Problems 2.\ref%
{prob2.1} and 2.\ref{prob2.2}) while the other factor is; we name them 
\textit{orient-} and \textit{skew-}factors, respectively. We show that a
given groupoid $\left( X,\bullet \right) $ with $x\bullet y\in \{x,y\}$, for
all $x,y$ in $X$, has a commutative method-2 factorization (Problem 3.\ref%
{prob3.2}). Section five briefly applies our two methods to some of the
algebras listed in section two; and discusses a promising relationship to
graph theory.

Finally, in our last section we generalize and summarize our findings that
certain groupoids/algebras decompose into distinct groupoids via $\left(
1\right) $ an operation on the parent groupoid and the left-zero-semigroup
simultaneously, which is a generalization of our first method; or $\left(
2\right) $ an operation which acts on the parent-groupoid and the
left-zero-semigroup separately, hence resulting in a generalization of our
second method.

Notions of \textquotedblleft method\textquotedblright -\textit{composite,}
\textquotedblleft method\textquotedblright -\textit{normal},
\textquotedblleft factor\textquotedblright -\textit{prime} and
\textquotedblleft partially\textquotedblright -\textit{left/right-prime}\
are used to classify and analyze various groupoids as well as other familiar
algebras. For simplicity, the left-zero-semigroup will be denoted as $%
id_{Bin\left( X\right) }$.

\section{Preliminaries}

\smallskip A \emph{groupoid} \cite{Br58} $\left( X,\bullet \right) $
consists of a non-empty set $X$ together with a binary operation $\bullet
:X\times X\rightarrow X$ where $x\bullet y\in X$ for all $x,y\in X$.

A groupoid $\left( X,\bullet \right) $ is \emph{strong }\cite{KN08} if and
only if for all $x,y\in X$, 
\begin{equation}
x\bullet y=y\bullet x\text{ implies }x=y.  \label{strong}
\end{equation}

A groupoid $\left( X,\bullet \right) $ is \emph{idempotent} if $x\bullet x=x$
for all $x\in X$.\medskip

\noindent \textbf{Example 2.1 }\cite{HKN10} Let $X=[0,\infty )$ and let $%
x\bullet y=\max \{0,x-y\}$ for any $x,y\in X$. Then $(X,\bullet )$ is a
strong groupoid. To visualize this, let's consider the associated Cayley
product table for \textquotedblleft $\bullet $\textquotedblright . For
simplicity, its partial table is displayed below which shows that $x\bullet
y=0$ for all $x\leq y$ and $x\bullet y\neq 0$ for all $x>y$:

\begin{equation*}
\begin{tabular}{c|cccccc}
$\bullet $ & 0 & 1 & 2 & 3 & 4 & $\cdots $ \\ \hline
\multicolumn{1}{c|}{0} & 0 & 0 & 0 & 0 & 0 & $\cdots $ \\ 
\multicolumn{1}{c|}{1} & 1 & 0 & 0 & 0 & 0 & $\cdots $ \\ 
\multicolumn{1}{c|}{2} & 2 & 1 & 0 & 0 & 0 & $\cdots $ \\ 
\multicolumn{1}{c|}{3} & 3 & 2 & 1 & 0 & 0 & $\cdots $ \\ 
\multicolumn{1}{c|}{4} & 4 & 3 & 2 & 1 & 0 & $\cdots $ \\ 
\multicolumn{1}{c|}{$\vdots $} & $\vdots $ & $\vdots $ & $\vdots $ & $\vdots 
$ & $\vdots $ & $\ddots $%
\end{tabular}%
\end{equation*}

Hence, the strong or anti-commutative property holds for all $x,y\in X$%
.\medskip

\noindent \textbf{Example 2.2 }\cite{HKN10} Let $X=\mathbf{%
\mathbb{R}
}$ be the set of all real numbers and let $x,y,e\in \mathbf{%
\mathbb{R}
}$. If we define a binary operation \textquotedblleft $\bullet $%
\textquotedblright\ on $X$ by $x\bullet y=(x-y)(x-e)+e$,

\noindent then the groupoid $(X,\bullet ,e)$ is not strong, since $%
x=e+\alpha ,$ $y=e-\alpha ,$ $\alpha \not=\pm e$ implies $x\bullet
y=y\bullet x$, but $x\not=y$.\bigskip

A groupoid $(X,\bullet )$ is a \emph{left-zero-semigroup} if $x\bullet y=x$
for all $x,y\in X$. Similarly, $(X,\bullet )$ is a \emph{right-zero-semigroup%
} if $x\bullet y=y$ for all $x,y\in X$. For the theory of semigroups, we
refer to \cite{CP61, NK96}. \medskip

A groupoid $(X,\bullet )$ is \emph{locally-zero} \cite{Fa11} if

(i) $x\bullet x=x$ for all $x\in X$; and

(ii) for any $x\not=y$ in $X$, $(\{x,y\},\bullet )$ is either a
left-zero-semigroup or a right-zero-semigroup.\bigskip

\noindent \textbf{Example 2.3} Given a set $X=\left\{ 0,1,2\right\} $, let
the binary operation \textquotedblleft $\bullet $\textquotedblright\ be
defined by the following Cayley product table:

\begin{equation*}
\begin{tabular}{l|lll}
$\bullet $ & 0 & 1 & 2 \\ \hline
0 & 0 & 0 & 2 \\ 
1 & 1 & 1 & 1 \\ 
2 & 0 & 2 & 2%
\end{tabular}%
\end{equation*}

\noindent Then the binary system $\left( X,\bullet \right) $ is locally-zero
and has the following subtables:%
\begin{equation*}
\text{%
\begin{tabular}{r|rr}
$\bullet $ & 0 & 1 \\ \hline
0 & 0 & 0 \\ 
1 & 1 & 1%
\end{tabular}%
\quad 
\begin{tabular}{r|rr}
$\bullet $ & 1 & 2 \\ \hline
1 & 1 & 1 \\ 
2 & 2 & 2%
\end{tabular}%
\quad 
\begin{tabular}{r|rr}
$\bullet $ & 0 & 2 \\ \hline
0 & 0 & 2 \\ 
2 & 0 & 2%
\end{tabular}%
}
\end{equation*}%
\noindent where $\left( \left\{ 0,1\right\} ,\bullet \right) $ is a
left-zero-semigroup; $\left( \left\{ 1,2\right\} ,\bullet \right) $ is also
a left-zero-semigroup; and $\left( \left\{ 0,2\right\} ,\bullet \right) $ is
a right-zero-semigroup.\bigskip

The notion of the semigroup $(Bin(X),\diamond )$ was introduced by J.
Neggers and H.S. Kim in \cite{KN08}. Given a non-empty set $X$, let $Bin(X)$
denote the collection of all groupoids $\left( X,\bullet \right) $, where $%
\bullet :X\times X\rightarrow X$ is a map. Given elements $(X,\ast )$ and $%
(X,\circ )$ of $Bin\left( X\right) $, define a binary product
\textquotedblleft $\diamond $\textquotedblright\ on these groupoids as
follows: 
\begin{equation}
(X,\ast )\diamond (X,\circ )=(X,\bullet )  \label{bin-grpd-prod}
\end{equation}

where 
\begin{equation}
x\bullet y=(x\ast y)\circ (y\ast x)  \label{bin-xy-prod}
\end{equation}%
for all $x$, $y\in X$. This turns $(Bin(X),\diamond )$ into a semigroup with
identity, the left-zero-semigroup, and an analog of negative one in the
right-zero-semigroup.

The present author \cite{Fa11} showed that a groupoid $\left( X,\bullet
\right) $ commutes, relative to the product \textquotedblleft $\diamond $%
\textquotedblright , if and only if any 2-element subset of $\left(
X,\bullet \right) $ is a subgroupoid that is either a left-zero-semigroup or
a right-zero-semigroup. Thus, $\left( X,\bullet \right) $ is an element of
the \emph{center} $ZBin(X)$ of the semigroup $(Bin(X),\diamond )$, defined
as follows:%
\begin{equation*}
ZBin\left( X\right) =\{\left( X,\bullet \right) \in Bin\left( X\right) \text{
}|\text{ }(X,\bullet )\,\diamond \,(X,\ast )\,=\,(X,\ast )\,\diamond
\,(X,\bullet ),\text{ }\forall (X,\ast )\in Bin(X)\}\text{.}
\end{equation*}%
In turn, several properties were obtained.\bigskip

\noindent \textbf{Theorem 2.4} \cite{KN08} \textsl{The collection }$%
(Bin(X),\diamond )$\textsl{\ of all binary systems (groupoids or algebras)
defined on }$X$\textsl{\ is a semigroup, i.e., the operation
\textquotedblleft }$\diamond $\textsl{\textquotedblright\ as defined in
general is associative. Furthermore, the left-zero-semigroup is an identity
for this operation.}\medskip

\noindent \textbf{Proposition 2.5} \cite{KN08}\textsl{\ Let }$(X,\bullet )$%
\textsl{\ be the right-zero-sermigroup on }$X$\textsl{. Then }$\left(
X,\bullet \right) \in $\textsl{\ }$Str(X)$\textsl{, the collection of all
strong groupoids on }$X$.\medskip

\noindent \textbf{Proposition 2.6} \cite{Fa11} \textsl{The left-zero
semigroup and right-zero semigroup on X are both in }$ZBin(X)$\textsl{.}%
\medskip

\noindent \textbf{Corollary 2.7.} \cite{Fa11} \textsl{The collection of all
locally-zero groupoids on }$X$\textsl{\ forms a subsemigroup of }$%
(Bin(X),\diamond )$\textsl{.}\medskip

\noindent \textbf{Proposition 2.8 }\cite{Fa11} \textsl{Let }$(X,\bullet )$%
\textsl{\ be a locally-zero groupoid. Then }$\left( X,\bullet \right)
\diamond \left( X,\bullet \right) =id_{Bin\left( X\right) }$,\textsl{\ the
left-zero-semigroup on }$X$\textsl{.}\medskip

Let $\left( X,\bullet \right) $ be an element of the semigroup $%
(Bin(X),\diamond )$, we say that $\left( X,\bullet \right) $ is a \textit{%
unit} if and only if there exists an element $\left( X,\ast \right) \in
Bin\left( X\right) $ such that 
\begin{equation}
(X,\bullet )\,\diamond \,(X,\ast )=id_{Bin\left( X\right) }=(X,\ast
)\,\diamond \,(X,\bullet )\text{.}  \label{unit}
\end{equation}

Subsequently, by Proposition 2.8, a locally-zero-groupoid is a unit in $%
Bin\left( X\right) $.\smallskip

The logic-based $BCK$/$BCI$-algebras were introduced by Is\'{e}ki and Imai
in \cite{IT78} as propositional calculus, but later in \cite{IK80} developed
into the present notion of $BCK$/$BCI$ which have since then been
investigated thoroughly by numerous researchers. J. Neggers and H. S. Kim
generalized a $BCK$-algebra \cite{MJ94} by introducing the notion of a $d$%
-algebra in \cite{NJK99}. They also introduced $B$-algebras in \cite{AKN02}.%
\textbf{\ }C. B. Kim and H. S. Kim generalized a $B$-algebra by defining a $%
BG$-algebra in \cite{KK08}.

An algebra $\left( X,\bullet ,0\right) $ of type $\left( 2,0\right) $ is a $%
B $\textit{-}\emph{algebra\ }\cite{AKN02}\textbf{\ }if for all $x,y,z\in X$,
it satisfies the following axioms:

\begin{description}
\item[B1] $x\bullet x=0$,

\item[B2] $x\bullet 0=x$, and

\item[B] $(x\bullet y)\bullet z=x\bullet \lbrack z\bullet (0\bullet y)]$.
\end{description}

An algebra $\left( X,\bullet ,0\right) $ of type $\left( 2,0\right) $ is a $%
BG$\textit{-}\emph{algebra\ }\cite{KK08}\textbf{\ }if for all $x,y,z\in X$,
it satisfies (B1), (B2), and

\begin{description}
\item[BG] $x=\left( x\bullet y\right) \bullet (0\bullet y)$.
\end{description}

An algebra $\left( X,\bullet ,0\right) $ of type $\left( 2,0\right) $ is a $%
BCI$\textit{-}\emph{algebra\ }\cite{Y06}\textbf{\ }if for all $x,y,z\in X$,
it satisfies (B2) and:

\begin{description}
\item[I] $\left( \left( x\bullet y\right) \bullet \left( x\bullet z\right)
\right) \bullet \left( z\bullet y\right) =0$,

\item[BH] $x\bullet y=0$ and $y\bullet x=0$ implies $x=y$. \smallskip
\end{description}

\noindent \textbf{Example 2.9} \cite{Y06} Let $X=\{0,1,a,b\}$. Define a
binary operation \textquotedblleft $\bullet $\textquotedblright\ on $X$ by
the following product table:%
\begin{equation*}
\begin{tabular}{c|cccc}
$\bullet $ & 0 & 1 & a & b \\ \hline
0 & 0 & 0 & a & a \\ 
1 & 1 & 0 & a & a \\ 
a & a & a & 0 & 0 \\ 
b & b & a & 1 & 0%
\end{tabular}%
\ 
\end{equation*}%
\smallskip

\noindent Then $(X,\bullet ,0)$ is a $BCI$-algebra.\bigskip

A $BCI$-algebra $\left( X,\bullet ,0\right) $ is a $BCK$-\emph{algebra} \cite%
{MJ94} if it satisfies the next additional axiom:\medskip

\begin{description}
\item[K] $0\bullet x=0$ for all $x\in X$.
\end{description}

An algebra $(X,\bullet ,0)$ of type $\left( 2,0\right) $ is a $d$\emph{%
-algebra} provided that for all $x,$ $y\in X$, it satisfies (B1), (K) and
(BH).\smallskip

A $d$-algebra is \emph{strong} if for all $x,$ $y$ $\in X$:

\begin{description}
\item[d-3$^{\prime }$] $x\bullet y=y\bullet x$ implies $x=y.$
\end{description}

Otherwise we consider the $d$-algebra to be \emph{exceptional}. For more
information on $d$-algebras we refer to \cite{AKN07, AKN11, NJK99, NK99}%
.\bigskip

\noindent \textbf{Example 2.10} \cite{NJK99} Let $(X,\bullet )=(%
\mathbb{Z}
_{5},\bullet )$ where \textquotedblleft $\bullet $\textquotedblright\ is
defined by the following Cayley table:%
\begin{equation*}
\begin{tabular}{c|ccccc}
$\bullet $ & 0 & 1 & 2 & 3 & 4 \\ \hline
0 & 0 & 0 & 0 & 0 & 0 \\ 
1 & 1 & 0 & 1 & 0 & 1 \\ 
2 & 2 & 2 & 0 & 3 & 0 \\ 
3 & 3 & 3 & 2 & 0 & 3 \\ 
4 & 4 & 4 & 1 & 1 & 0%
\end{tabular}%
\end{equation*}%
\smallskip

\noindent Then $(%
\mathbb{Z}
_{5},\bullet ,0)$ is a $d$-algebra which is not a $BCK$-algebra. For details
on $BCK$-algebras, see \cite{Ior08, MJ94, Y06}.\bigskip

Y. B. Jun, E. H. Roh and H. S. Kim in \cite{JRK98} introduced the notion of
a $BH$-algebra which is a generalization of $BCK/BCI/BCH$-algebras. There
are many other generalizations of similar algebras. We summarize several
properties which are used as axioms to define each algebraic structure . Let 
$(X,\bullet ,0)$ be an algebra of type $\left( 2,0\right) $, for any $x,$ $%
y, $ $z\in X$:

\begin{description}
\item[B1] $x\bullet x=0$,

\item[B2] $x\bullet 0=x$,

\item[B] $(x\bullet y)\bullet z=x\bullet (z\bullet (0\bullet y)),$

\item[BG] $x=(x\bullet y)\bullet (0\bullet y),$

\item[BM] $(z\bullet x)\bullet (z\bullet y)=y\bullet x,$

\item[BH] $x\bullet y=0$ and $y\bullet x=0\Rightarrow x=y$,

\item[BF] $0\bullet (x\bullet y)=y\bullet x$,

\item[BN] $(x\bullet y)\bullet z=(0\bullet z)\bullet (y\bullet x)$,

\item[BO] $x\bullet (y\bullet z)=(x\bullet y)\bullet (0\bullet z)$,

\item[BP1] $x\bullet (x\bullet y)=y$,

\item[BP2] $(x\bullet z)\bullet (y\bullet z)=x\bullet y$,

\item[Q] $(x\bullet y)\bullet z=(x\bullet z)\bullet y$,

\item[CO] $(x\bullet y)\bullet z=x\bullet (y\bullet z)$,

\item[BZ] $((x\bullet z)\bullet (y\bullet z))\bullet (x\bullet y)=0$,

\item[K] $0\bullet x=0$,

\item[I] $((x\bullet y)\bullet (x\bullet z))\bullet (z\bullet y)=0$,

\item[BI] $x\bullet (y\bullet x)=x$.\medskip
\end{description}

An algebra $(X,\bullet ,0)$ of type $(2,0)$ is classified according to a
combination of the above axioms as noted in \textquotedblleft \textbf{Figure
1}\textquotedblright\ below. For instance, $(X,\bullet ,0)$ is a $BI$%
-algebra \cite{SKR17} if satisfies in (B1) and (BI). For detailed
information on each, please see [2-6, 14-26, 31, 32, 34, 36].

\begin{figure}[t]
\centering
\includegraphics[width=1\linewidth, height=0.5\textheight]{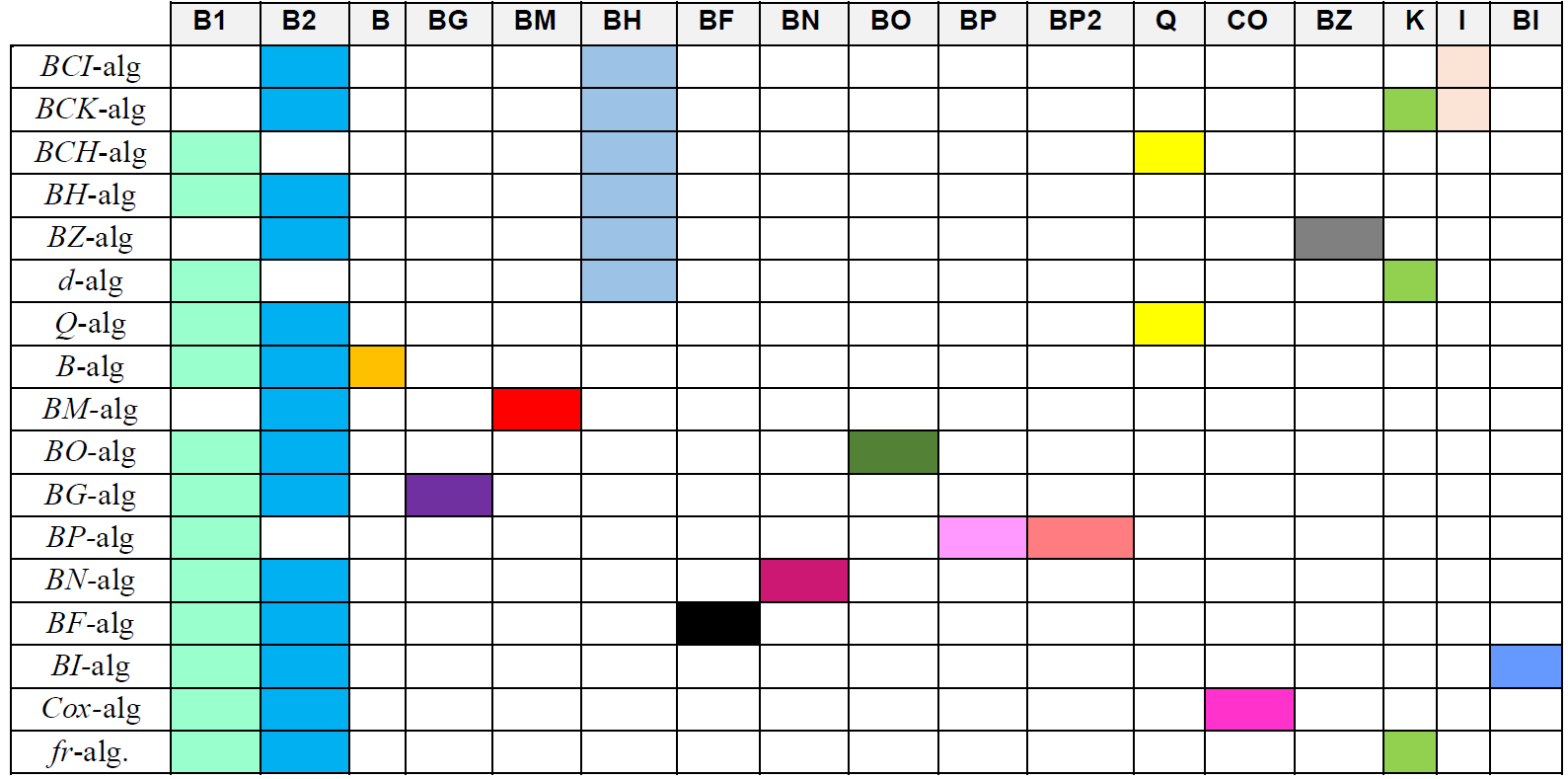}
\caption{Comparison of Algebras}
\label{fig:tableofalgebras}
\end{figure}
\newpage

\section{Similar-Signature Factorization}

In this section, we present a unique factorization of a given groupoid by
\textquotedblleft deriving\textquotedblright\ two factors from it and from
the left-zero-semigroup simultaneously.\smallskip

Let $(X,\bullet )$ be a groupoid of finite order, i.e., $\left\vert
X\right\vert =n$. Then $d^{\bullet }$ is the \emph{diagonal function} of $%
(X,\bullet )\ $such that $d^{\bullet }:%
\mathbb{N}
\longrightarrow X$ where $d^{\bullet }(i)=x_{i}\bullet x_{i}$, $i=1,2,...,n$
for all $x_{i}$ $\in X$.\smallskip

\noindent \textbf{Example 3.1 }Let $\left( X,\bullet ,0\right) $ and $\left(
X,\ast \right) $ be a $d$-algebra and an idempotent algebra, respectively.
Then $x\bullet x=0$ and $x\ast x=x$; or $d^{\bullet }=0$ and $d^{\ast }=x$
for all $x\in X$. \medskip

Two binary systems $(X,\ast )$ and $(X,\bullet )$ are said to be \emph{%
similar} if they have the same diagonal function, that is, $d^{\ast
}=d^{\bullet }$.\smallskip

Two binary systems $(X,\ast )$ and $(X,\bullet )$ are said to be \emph{%
signature }if

(i) $x\ast y=x\bullet y$ when $x\neq y;$ and

(ii) $x\ast x\neq x\bullet x$ for all $x\in X$.\medskip

Let $\left( X,\bullet \right) $ be a groupoid. \textit{Derive} groupoids $%
(X,\ast )$ and $(X,\circ )$ from $\left( X,\bullet \right) $ and $%
id_{Bin\left( X\right) }$, simultaneously, such that for all $x,y\in X$,\ 
\begin{equation}
x\ast y=%
\begin{cases}
x & \text{if $x=y$,} \\ 
x\bullet y & \text{otherwise.}%
\end{cases}%
\text{\quad and\quad }x\circ y=%
\begin{cases}
x\bullet x & \text{if $x=y$,} \\ 
x & \text{otherwise.}%
\end{cases}
\label{Sig-Sim-F's}
\end{equation}%
\ \ The groupoids $\left( X,\ast \right) $ and $\left( X,\circ \right) $ are
said to be the \emph{signature- }and the \emph{similar-factors }of $\left(
X,\bullet \right) $, respectively, denoted by $U\left( X,\bullet \right) $
and $A(X,\bullet )$. The product \textquotedblleft $\diamond $%
\textquotedblright\ is associative but not commutative. Hence, for $\left(
X,\bullet \right) \in Bin\left( X\right) $, we may have a \emph{$UA$%
-factorization} such that%
\begin{equation}
\left( X,\bullet \right) =U\left( X,\bullet \right) \diamond A\left(
X,\bullet \right)  \label{UA-prod}
\end{equation}%
or an \emph{$AU$-factorization} such that%
\begin{equation}
\left( X,\bullet \right) =A\left( X,\bullet \right) \diamond U\left(
X,\bullet \right) \text{.}  \label{AU-prod}
\end{equation}

By the equations in \ref{Sig-Sim-F's}, it follows that for any given
groupoid $\left( X,\bullet \right) $,

\begin{enumerate}
\item $U\left( X,\bullet \right) $ is similar to $id_{Bin\left( X\right) }$
while $A\left( X,\bullet \right) $ is similar to $\left( X,\bullet \right) $%
; and

\item $U\left( X,\bullet \right) $ is signature with $\left( X,\bullet
\right) $ while $A\left( X,\bullet \right) $ is signature with $%
id_{Bin\left( X\right) }$.\medskip
\end{enumerate}

\noindent \textbf{Proposition 3.2} \textsl{The similar-factor of a groupoid
is strong.}

\noindent \textit{Proof.} Given $\left( X,\bullet \right) $ $\in Bin\left(
X\right) $, let $\left( X,\circ \right) =A\left( X,\bullet \right) $.

\begin{enumerate}
\item[(i)] If $x=y$, then $x\circ y=x\circ x=x\bullet x=y\bullet y=y\circ
y=y\circ x$.

\item[(ii)] If $x\neq y$ and $x\circ y=y\circ x$ for any $x,y\in X$. Then $%
x\circ y=x$ and $y\circ x=y$. Thus, $x=y$, a contradiction.
\end{enumerate}

\noindent Therefore, $\left( X,\circ \right) $ is strong.

\begin{flushright}
$\blacksquare $
\end{flushright}

\noindent \textbf{Example 3.3 } Let $\left( X,\bullet ,0\right) $ be the $%
BCI $-algebra defined in Example 2.9. In accordance with equation \ref%
{Sig-Sim-F's}, derive its \textit{signature-} and \textit{similar-} factors $%
U(X,\bullet ,0)$ and $A(X,\bullet ,0)$, respectively. Let groupoids $\left(
X,\ast ,0\right) :=U(X,\bullet ,0)$ and $\left( X,\circ ,0\right) :=A\left(
X,\bullet ,0\right) $ be given. We obtain:%
\begin{equation*}
\begin{tabular}{c|cccc}
$\ast $ & 0 & 1 & a & b \\ \hline
0 & 0 & 0 & a & a \\ 
1 & 1 & 1 & a & a \\ 
a & a & a & a & 0 \\ 
b & b & a & 1 & b%
\end{tabular}%
\quad \text{and}\quad 
\begin{tabular}{c|cccc}
$\circ $ & 0 & 1 & a & b \\ \hline
0 & 0 & 0 & 0 & 0 \\ 
1 & 1 & 0 & 1 & 1 \\ 
a & a & a & 0 & a \\ 
b & b & b & b & 0%
\end{tabular}%
\end{equation*}

\noindent It remains to verify that $\left( X,\bullet ,0\right) =\left(
X,\ast ,0\right) \diamond \left( X,\circ ,0\right) $ and/or $\left(
X,\bullet ,0\right) =\left( X,\circ ,0\right) \diamond \left( X,\ast
,0\right) $. This will be discussed in more detail in the next section.
However, there is a very interesting fact in this example: the two factors
are distinct from each other, their parent groupoid, and the
left-zero-semigroup. In summary:

\begin{enumerate}
\item $\left( X,\ast ,0\right) \neq \left( X,\circ ,0\right) $; (Problem 1.%
\ref{prob1.1})

\item $\left( X,\ast ,0\right) \neq \left( X,\bullet ,0\right) \neq \left(
X,\circ ,0\right) $; (Problem 1.\ref{prob1.3})

\item $\left( X,\ast ,0\right) \neq id_{Bin\left( X\right) }\neq \left(
X,\circ ,0\right) $. (Problem 1.\ref{prob1.4})
\end{enumerate}

\noindent This is important since it is not always the case that all three
distinctions hold as the following example demonstrates.\medskip

\noindent \textbf{Example 3.4 }Let $(X,\bullet )=(%
\mathbb{Z}
_{3},\bullet )$ where \textquotedblleft $\bullet $\textquotedblright\ is
defined by the following Cayley table:%
\begin{equation*}
\begin{tabular}{c|ccc}
$\bullet $ & 0 & 1 & 2 \\ \hline
0 & 0 & 0 & 0 \\ 
1 & 1 & 0 & 1 \\ 
2 & 2 & 2 & 0%
\end{tabular}%
\ 
\end{equation*}

\noindent Then $\left( X,\bullet ,0\right) $ is a $BI$-algebra. Derive its 
\textit{signature}- and \textit{similar-}factors $U(X,\bullet ,0)$ and $%
A(X,\bullet ,0)$, respectively, in accordance to the equations in \ref%
{Sig-Sim-F's}. Let $\left( X,\ast ,0\right) :=U\left( X,\bullet ,0\right) $
and $\left( X,\circ ,0\right) :=A\left( X,\bullet ,0\right) $, hence:%
\begin{equation*}
\begin{tabular}{c|ccc}
$\ast $ & 0 & 1 & 2 \\ \hline
0 & 0 & 0 & 0 \\ 
1 & 1 & 1 & 1 \\ 
2 & 2 & 2 & 2%
\end{tabular}%
\quad \text{and}\quad 
\begin{tabular}{c|ccc}
$\circ $ & 0 & 1 & 2 \\ \hline
0 & 0 & 0 & 0 \\ 
1 & 1 & 0 & 1 \\ 
2 & 2 & 2 & 0%
\end{tabular}%
\end{equation*}

\noindent Here we observe immediately that the \textit{similar}-factor $%
\left( X,\circ ,0\right) $ is equal to $\left( X,\bullet ,0\right) $ and the 
\textit{signature-}factor $\left( X,\ast ,0\right) $ is equal to $%
id_{Bin\left( X\right) }$. Thus this decomposition is basically a trivial
factorization, i.e., 
\begin{equation*}
\left( X,\bullet ,0\right) =\left( X,\ast ,0\right) \diamond \left( X,\circ
,0\right) =id_{Bin\left( X\right) }\diamond \left( X,\bullet ,0\right)
\end{equation*}%
\noindent and

\begin{equation*}
\left( X,\bullet ,0\right) =\left( X,\circ ,0\right) \diamond \left( X,\ast
,0\right) =\left( X,\bullet ,0\right) \diamond id_{Bin\left( X\right) }\text{%
.}
\end{equation*}

\subsection{$UA$-Factorization}

In this subsection, we explore a $UA$-factorization of a given groupoid $%
\left( X,\bullet \right) $ in $Bin\left( X\right) $. In the next subsection,
a $AU$-factorization is considered, where the order of the product of the
two factors is \textquotedblleft reversed\textquotedblright . We emphasize
that such factorization is unique and not necessarily reversible. Then, we
classify a given groupoid as $UA$\textit{- }and/or\textit{\ }$AU$-\textit{%
composite}, $u$-\textit{composite} or\textit{\ }$u$\textit{-normal; }and as%
\textit{\ signature- }or\textit{\ similar-prime. }\medskip

\noindent \textbf{Example 3.1.1} Let $X=%
\mathbb{Z}
$ be the set of all integers and let \textquotedblleft $-$%
\textquotedblright\ be the usual subtraction on $%
\mathbb{Z}
$. Then $\left( 
\mathbb{Z}
,-\right) $ is a $BH$-algebra since it satisfies axioms \textbf{B1}, \textbf{%
B2} and \textbf{BH} as seen from its partial table below: 
\begin{equation*}
\begin{tabular}{c|ccccccccc}
$-$ & $\cdots $ & -2 & -1 & 0 & 1 & 2 & 3 & 4 & $\cdots $ \\ \hline
$\vdots $ & $\ddots $ & $\vdots $ & $\vdots $ & $\vdots $ & $\vdots $ & $%
\vdots $ & $\vdots $ & $\vdots $ & $\cdots $ \\ 
-2 & 1 & 0 & -1 & -2 & -3 & -4 & -5 & -6 & $\cdots $ \\ 
-1 & 2 & 1 & 0 & -1 & -2 & -3 & -4 & -5 & $\cdots $ \\ 
0 & 3 & 2 & 1 & 0 & -1 & -2 & -3 & -4 & $\cdots $ \\ 
1 & 4 & 3 & 2 & 1 & 0 & -1 & -2 & -3 & $\cdots $ \\ 
2 & 5 & 4 & 3 & 2 & 1 & 0 & -1 & -2 & $\cdots $ \\ 
3 & 6 & 5 & 4 & 3 & 2 & 1 & 0 & -1 & $\cdots $ \\ 
4 & 7 & 6 & 5 & 4 & 3 & 2 & 1 & 0 & $\cdots $ \\ 
$\vdots $ & $\vdots $ & $\vdots $ & $\vdots $ & $\vdots $ & $\vdots $ & $%
\vdots $ & $\vdots $ & $\vdots $ & $\ddots $%
\end{tabular}%
\end{equation*}

\noindent Define two binary operations \textquotedblleft $\ast $%
\textquotedblright\ and \textquotedblleft $\circ $\textquotedblright\ on $%
\mathbb{Z}
$ such that for all $x,y\in 
\mathbb{Z}
$, 
\begin{equation*}
x\ast y=%
\begin{cases}
x & \text{if $x=y$,} \\ 
x-y & \text{otherwise.}%
\end{cases}%
\text{\quad and}\quad x\circ y=%
\begin{cases}
0 & \text{if $x=y$,} \\ 
x & \text{otherwise.}%
\end{cases}%
\noindent
\end{equation*}

\noindent Then it is easy to check that $\left( 
\mathbb{Z}
,-\right) =(%
\mathbb{Z}
,\ast )\diamond (%
\mathbb{Z}
,\circ )$ and $(%
\mathbb{Z}
,\ast )=U(%
\mathbb{Z}
,-)$ and $(%
\mathbb{Z}
,\circ )=A(%
\mathbb{Z}
,-)$. Thus we have a $UA$-factorization of $\left( 
\mathbb{Z}
,-\right) .\medskip $

A groupoid $\left( X,\bullet \right) $ is said to be \emph{signature-prime}%
\textit{\ }if $U\left( X,\bullet \right) $ $=id_{Bin\left( X\right) }$, and
is said to be \emph{similar-prime} if $A\left( X,\bullet \right) $ $%
=id_{Bin\left( X\right) }$. Alternatively, if $\left( X,\bullet \right) $ is
neither \textit{signature}- nor \textit{similar}-prime, then $\left(
X,\bullet \right) $ is said to be

\qquad (1) \emph{UA-composite} if $\left( X,\bullet \right) =U\left(
X,\bullet \right) \diamond A\left( X,\bullet \right) $;

\qquad (2) \emph{$AU$-composite} if $\left( X,\bullet \right) =A\left(
X,\bullet \right) \diamond U\left( X,\bullet \right) $.

Consequently, $\left( X,\bullet \right) $ is said to be $u$\emph{-composite}
if both $\left( 1\right) $ and $\left( 2\right) $ hold.\medskip

\noindent \textbf{Example 3.1.2} Let $(X,\bullet )=(%
\mathbb{Z}
_{5},\bullet )$ where the product \textquotedblleft $\bullet $%
\textquotedblright\ is defined by the following Cayley table: 
\begin{equation*}
\begin{tabular}{c|ccccc}
$\bullet $ & 0 & 1 & 2 & 3 & 4 \\ \hline
0 & 3 & 2 & 2 & 1 & 1 \\ 
1 & 1 & 3 & 3 & 2 & 3 \\ 
2 & 3 & 3 & 0 & 3 & 0 \\ 
3 & 1 & 0 & 1 & 1 & 2 \\ 
4 & 1 & 1 & 2 & 4 & 2%
\end{tabular}%
\ 
\end{equation*}

\noindent If we derive its \textit{signature}- and \textit{similar-} factors 
$\left( 
\mathbb{Z}
_{5},\ast \right) =U(%
\mathbb{Z}
_{5},\bullet )$ and $A(%
\mathbb{Z}
_{5},\bullet )=\left( 
\mathbb{Z}
_{5},\circ \right) $ as in (\ref{Sig-Sim-F's}), then we have their $\diamond 
$ product as follows:%
\begin{equation*}
\text{\quad }%
\begin{tabular}{c|ccccc}
$\ast $ & 0 & 1 & 2 & 3 & 4 \\ \hline
0 & 0 & 2 & 2 & 1 & 1 \\ 
1 & 1 & 1 & 3 & 2 & 3 \\ 
2 & 3 & 3 & 2 & 3 & 0 \\ 
3 & 1 & 0 & 1 & 3 & 2 \\ 
4 & 1 & 1 & 2 & 4 & 4%
\end{tabular}%
\quad \diamond \quad 
\begin{tabular}{c|ccccc}
$\circ $ & 0 & 1 & 2 & 3 & 4 \\ \hline
0 & 3 & 0 & 0 & 0 & 0 \\ 
1 & 1 & 3 & 1 & 1 & 1 \\ 
2 & 2 & 2 & 0 & 2 & 2 \\ 
3 & 3 & 3 & 3 & 1 & 3 \\ 
4 & 4 & 4 & 4 & 4 & 2%
\end{tabular}%
\quad =\quad 
\begin{tabular}{c|ccccc}
$\nabla $ & 0 & 1 & 2 & 3 & 4 \\ \hline
0 & 3 & 2 & 2 & 3 & 3 \\ 
1 & 1 & 3 & 1 & 2 & 3 \\ 
2 & 3 & 1 & 0 & 3 & 0 \\ 
3 & 3 & 0 & 1 & 1 & 2 \\ 
4 & 3 & 1 & 2 & 4 & 2%
\end{tabular}%
\quad
\end{equation*}

\noindent We can clearly conclude that $U(%
\mathbb{Z}
_{5},\bullet )\diamond A(%
\mathbb{Z}
_{5},\bullet )\neq \left( 
\mathbb{Z}
_{5},\bullet \right) $ since $\left( 
\mathbb{Z}
_{5},\bullet \right) \neq \left( 
\mathbb{Z}
_{5},\nabla \right) $ and hence such a groupoid does not have a $UA$%
-factorization. Moreover, $(%
\mathbb{Z}
_{5},\bullet )$ is not a strong groupoid since $0\bullet 4=4\bullet 0$. In
turn, we have the next theorem.\medskip

\noindent \textbf{Theorem 3.1.3} \textsl{A strong groupoid has a }$UA$%
\textsl{-factorization.}\medskip

\noindent \textit{Proof.} Let $(X,\bullet )$ $\in Str\left( X\right) $, the
collection of all strong groupoids defined on $X$, and let $(X,\odot
)=(X,\ast )$ $\diamond (X,\circ )$ where $(X,\ast )=U(X,\bullet )$ and $%
(X,\circ )=A(X,\bullet )$. Then $x\odot y=(x\ast y)\circ (y\ast x)$ for all $%
x,y\in X$. It follows that $x\ast x=x,$ $x\ast y=x\bullet y$ when $x\not=y$;
and $x\circ x=x\bullet x$, $x\circ y=x$ when $x\not=y$.

\noindent Next, we show that $(X,\bullet )=(X,\odot )$. Given $x,y\in X$, if 
$x=y$, then $x\odot x=(x\ast x)\circ (x\ast x)=x\circ x=x\bullet x$. Assume $%
x\not=y$, we claim that $x\ast y=y\ast x$ is not possible:

(i) If $x\ast y=y\ast x$, then $x\bullet y=x\ast y=y\ast x=y\bullet x$.
Since $(X,\bullet )$ is strong, we obtain $x=y$, a contradiction.

(ii) If $x\ast y\not=y\ast x$, then $x\ast y=x\bullet y,$ $y\ast x=y\bullet
x $, since $x\not=y$.

\noindent Therefore $x\odot y=(x\ast y)\circ (y\ast x)=(x\bullet y)\circ
(y\bullet x)=x\bullet y$, since $x\bullet y\not=y\bullet x$. This proves
that $(X,\odot )=(X,\bullet )$.

\begin{flushright}
$\blacksquare $
\end{flushright}

\noindent \textbf{Corollary 3.1.4} \textsl{The factorization in Theorem
3.1.3 is unique.} \medskip

\noindent \textit{Proof.} Let $(X,\bullet )$ be a strong groupoid with a $UA$%
-factorization such that $(X,\bullet )=(X,\ast )$ $\diamond (X,\circ )$
where $(X,\ast )=U(X,\bullet )$ and $(X,\circ )=A(X,\bullet )$. Let $%
(X,\bullet )=(X,\bigtriangledown )$ $\diamond (X,\bigtriangleup )$ where $%
(X,\bigtriangledown )=U(X,\bullet )$ and $(X,\bigtriangleup )=A(X,\bullet )$%
. For any $x\in X$, we have $x\ast x=x=x\bigtriangledown x$, and $x\ast
y=x\bigtriangledown y$ when $x\not=y$. Hence $(X,\ast )=(X,\bigtriangledown
) $. Similarly, if $x\in X$, then $x\circ x=x\bullet x=x\bigtriangleup x$.
When $x\not=y$, we have $x\circ y=x=x\bigtriangleup y$, proving that $%
(X,\circ )=(X,\bigtriangleup )$.

\begin{flushright}
$\blacksquare $
\end{flushright}

\noindent \textbf{Example 3.1.5 }\cite{NJK99} Consider the $d$-algebra $%
\left( X,\bullet ,0\right) $ from Example 2.10. Observe that $\left(
X,\bullet ,0\right) $ is a strong $d$-algebra. Let $\left( X,\ast ,0\right)
:=U(X,\bullet ,0)$ and $\left( X,\circ ,0\right) :=A\left( X,\bullet
,0\right) $, such that $U(X,\bullet ,0)$ and $A(X,\bullet ,0)$ are its
derived \textit{signature-} and \textit{similar-}factors, respectively, as
in $\left( \text{\ref{Sig-Sim-F's}}\right) $. Next, verify that $\left(
X,\ast ,0\right) \diamond \left( X,\circ ,0\right) =$\ $\left( X,\bullet
,0\right) $: 
\begin{equation*}
\text{\quad }%
\begin{tabular}{c|ccccc}
$\ast $ & 0 & 1 & 2 & 3 & 4 \\ \hline
0 & 0 & 0 & 0 & 0 & 0 \\ 
1 & 1 & 1 & 1 & 0 & 1 \\ 
2 & 2 & 2 & 2 & 3 & 0 \\ 
3 & 3 & 3 & 2 & 3 & 3 \\ 
4 & 4 & 4 & 1 & 1 & 4%
\end{tabular}%
\text{\quad }\diamond \text{\quad }%
\begin{tabular}{c|ccccc}
$\circ $ & 0 & 1 & 2 & 3 & 4 \\ \hline
0 & 0 & 0 & 0 & 0 & 0 \\ 
1 & 1 & 0 & 1 & 1 & 1 \\ 
2 & 2 & 2 & 0 & 2 & 2 \\ 
3 & 3 & 3 & 3 & 0 & 3 \\ 
4 & 4 & 4 & 4 & 4 & 0%
\end{tabular}%
\text{\quad }=\text{\quad }%
\begin{tabular}{c|ccccc}
$\bullet $ & 0 & 1 & 2 & 3 & 4 \\ \hline
0 & 0 & 0 & 0 & 0 & 0 \\ 
1 & 1 & 0 & 1 & 0 & 1 \\ 
2 & 2 & 2 & 0 & 3 & 0 \\ 
3 & 3 & 3 & 2 & 0 & 3 \\ 
4 & 4 & 4 & 1 & 1 & 0%
\end{tabular}%
\end{equation*}%
\noindent \noindent Indeed we can see that $x\bullet y=\left( x\ast y\right)
\circ \left( y\ast x\right) $ for any $x,y\in X$. For instance:%
\begin{eqnarray*}
\left( 1\ast 0\right) \circ \left( 0\ast 1\right) &=&1\circ 0=1=1\bullet 0%
\text{,} \\
\left( 3\ast 4\right) \circ \left( 4\ast 3\right) &=&3\circ 4=3=3\bullet 4%
\text{.}
\end{eqnarray*}%
\noindent Moreover, since $U\left( X,\bullet ,0\right) \neq id_{Bin\left(
X\right) }$ and $A\left( X,\bullet ,0\right) \neq id_{Bin\left( X\right) }$,
then $\left( X,\bullet ,0\right) $ is $UA$-composite.

\subsection{$AU$-Factorization}

In this subsection we reverse the order of the \textit{signature}- and 
\textit{similar}-factors of any groupoid $\left( X,\bullet \right) $ in $%
Bin\left( X\right) $. We conclude that an arbitrary groupoid $\left(
X,\bullet \right) $ will always have an $AU$-factorization. However, this
factorization might be trivial and hence the groupoid is either noted as 
\textit{signature-} or \textit{similar-}prime. Otherwise, if the
decomposition is not trivial, we say the groupoid is $AU$-composite.\bigskip

\noindent \textbf{Example 3.2.1} Let $\left( X,\bullet ,0\right) $ be the
strong $d$-algebra defined in Examples 2.10 and 3.1.5 in which we determined
that $\left( X,\bullet ,0\right) $ is $UA$-composite. Similarly, we can take
the product of $A\left( X,\bullet ,0\right) $ and $U\left( X,\bullet
,0\right) $ as follows:%
\begin{equation*}
\text{\quad }%
\begin{tabular}{c|ccccc}
$\circ $ & 0 & 1 & 2 & 3 & 4 \\ \hline
0 & 0 & 0 & 0 & 0 & 0 \\ 
1 & 1 & 0 & 1 & 1 & 1 \\ 
2 & 2 & 2 & 0 & 2 & 2 \\ 
3 & 3 & 3 & 3 & 0 & 3 \\ 
4 & 4 & 4 & 4 & 4 & 0%
\end{tabular}%
\text{\quad }\diamond \text{\quad }%
\begin{tabular}{c|ccccc}
$\ast $ & 0 & 1 & 2 & 3 & 4 \\ \hline
0 & 0 & 0 & 0 & 0 & 0 \\ 
1 & 1 & 1 & 1 & 0 & 1 \\ 
2 & 2 & 2 & 2 & 3 & 0 \\ 
3 & 3 & 3 & 2 & 3 & 3 \\ 
4 & 4 & 4 & 1 & 1 & 4%
\end{tabular}%
\text{\quad }=\text{\quad }%
\begin{tabular}{c|ccccc}
$\bullet $ & 0 & 1 & 2 & 3 & 4 \\ \hline
0 & 0 & 0 & 0 & 0 & 0 \\ 
1 & 1 & 0 & 1 & 0 & 1 \\ 
2 & 2 & 2 & 0 & 3 & 0 \\ 
3 & 3 & 3 & 2 & 0 & 3 \\ 
4 & 4 & 4 & 1 & 1 & 0%
\end{tabular}%
\end{equation*}%
\noindent By routine checking of $\left( x\circ y\right) \ast \left( y\circ
x\right) =x\bullet y$ for any $x,y\in X$, we conclude that $\left( X,\bullet
,0\right) $ has an $AU$-factorization. Moreover, we can see that this
particular groupoid has both, a non-trivial $UA$- and $AU$-factorization.
Therefore, $\left( X,\bullet ,0\right) $ is $u$-composite.\medskip

\noindent \textbf{Remark 3.2.2} Note that $A(X,\bullet )\diamond U(X,\bullet
)=U(X,\bullet )\diamond A(X,\bullet )$ does not imply that $\left( X,\bullet
\right) $ is $u$-composite. It simply implies that the factors of $\left(
X,\bullet \right) $ commute. This motivates the next definition.\medskip

A groupoid $\left( X,\bullet \right) $ is said to be $u$\emph{-normal}%
\textit{\ }if it admits a $UA$- and an $AU$-factorization, i.e., if

\qquad (i) $\left( X,\bullet \right) =U\left( X,\bullet \right) \diamond
A\left( X,\bullet \right) $, and

\qquad (ii) $\left( X,\bullet \right) =A\left( X,\bullet \right) \diamond
U\left( X,\bullet \right) $.\medskip

\noindent \textbf{Theorem 3.2.3} \textsl{Any given groupoid has an }$AU$%
\textsl{-factorization, i.e., if }$(X,\bullet )\in Bin(X)$\textsl{$,$ then}%
\begin{equation*}
\left( X,\bullet \right) =A(X,\bullet )\diamond U(X,\bullet )\text{.}
\end{equation*}%
\smallskip \noindent \textit{Proof.} Let $(X,\bullet )$ $\in Bin\left(
X\right) $ and let $(X,\odot )=(X,\circ )$ $\diamond (X,\ast )$ where $%
(X,\ast )=U(X,\bullet )$ and $(X,\circ )=A(X,\bullet )$. Then $x\odot
y=(x\circ y)\ast (y\circ x)$ for all $x,y\in X$. It follows that $x\ast x=x,$
$x\ast y=x\bullet y$ when $x\not=y$, and $x\circ x=x\bullet x,$ $x\circ y=x$
when $x\not=y$. Given $x,y\in X$, if $x=y$, then $x\odot x=(x\circ x)\ast
(x\circ x)=\left( x\bullet x\right) \ast \left( x\bullet x\right) =x\bullet
x $. Assume $x\not=y$, then $x\odot y=(x\circ y)\ast (y\circ x)=x\ast
y=x\bullet y$. This proves that $(X,\odot )=(X,\bullet )$.

\begin{flushright}
$\blacksquare $
\end{flushright}

\noindent \textbf{Corollary 3.2.4} \textsl{The factorization in Theorem
3.2.3 is unique.\smallskip }

\noindent \textit{Proof.} The proof is similar to that of Corollary 3.1.4.

\begin{flushright}
$\blacksquare $
\end{flushright}

\noindent \textbf{Corollary 3.2.5 }\textsl{A strong groupoid is }$u$\textsl{%
-normal.\smallskip }

\noindent \textit{Proof.} The proof follows directly from Theorems 3.1.3,
3.2.3 and the definition.

\begin{flushright}
$\blacksquare $
\end{flushright}

\noindent \textbf{Example 3.2.6} Let $\left( X,\bullet \right) =\left(
\left\{ 0,1,2\right\} ,+\right) $ be the cyclic group of order 3. Observe
that $\left( \left\{ 0,1,2\right\} ,+\right) $ has an $AU$-factorization but
fails to have a $UA$-factorization. Take $\left( \left\{ 0,1,2\right\} ,\ast
\right) =$ $U\left( \left\{ 0,1,2\right\} ,+\right) $ and $\left( \left\{
0,1,2\right\} ,\circ \right) =A\left( \left\{ 0,1,2\right\} ,+\right) $ such
that:%
\begin{equation*}
x\circ y=%
\begin{cases}
(x+x)\text{ mod }3 & \text{ if }x=y, \\ 
x & \text{otherwise.}%
\end{cases}%
\text{\quad and\quad }\ x\ast y=%
\begin{cases}
x; & \text{ if }x=y, \\ 
(x+y)\text{ mod }3 & \text{otherwise.}%
\end{cases}%
\end{equation*}%
\noindent Routine checking of the product $A\left( \left\{ 0,1,2\right\}
,+\right) \diamond U\left( \left\{ 0,1,2\right\} ,+\right) $ gives $\left(
\left\{ 0,1,2\right\} ,+\right) $:%
\begin{equation*}
\text{\quad }%
\begin{tabular}{r|rrr}
$\circ $ & $0$ & $1$ & $2$ \\ \hline
$0$ & $0$ & $0$ & $0$ \\ 
$1$ & $1$ & $2$ & $1$ \\ 
$2$ & $2$ & $2$ & $1$%
\end{tabular}%
\text{\quad }\diamond \text{\quad }%
\begin{tabular}{r|rrr}
$\ast $ & $0$ & $1$ & $2$ \\ \hline
$0$ & $0$ & $1$ & $2$ \\ 
$1$ & $1$ & $1$ & $0$ \\ 
$2$ & $2$ & $0$ & $2$%
\end{tabular}%
\text{\quad }=\text{\quad }%
\begin{tabular}{r|rrr}
$+$ & $0$ & $1$ & $2$ \\ \hline
$0$ & $0$ & $1$ & $2$ \\ 
$1$ & $1$ & $2$ & $0$ \\ 
$2$ & $2$ & $0$ & $1$%
\end{tabular}%
\end{equation*}%
\noindent But, the product $U\left( \left\{ 0,1,2\right\} ,+\right) \diamond
A\left( \left\{ 0,1,2\right\} ,+\right) $ does not give $\left( \left\{
0,1,2\right\} ,+\right) $:%
\begin{equation*}
\text{\quad }%
\begin{tabular}{r|rrr}
$\ast $ & $0$ & $1$ & $2$ \\ \hline
$0$ & $0$ & $1$ & $2$ \\ 
$1$ & $1$ & $1$ & $0$ \\ 
$2$ & $2$ & $0$ & $2$%
\end{tabular}%
\text{\quad }\diamond \text{\quad }%
\begin{tabular}{r|rrr}
$\circ $ & $0$ & $1$ & $2$ \\ \hline
$0$ & $0$ & $0$ & $0$ \\ 
$1$ & $1$ & $2$ & $1$ \\ 
$2$ & $2$ & $2$ & $1$%
\end{tabular}%
\text{\quad }=\text{\quad }%
\begin{tabular}{r|rrr}
$\nabla $ & $0$ & $1$ & $2$ \\ \hline
$0$ & $0$ & $2$ & $1$ \\ 
$1$ & $2$ & $2$ & $0$ \\ 
$2$ & $1$ & $0$ & $1$%
\end{tabular}%
\end{equation*}%
\noindent Therefore, $\left( \left\{ 0,1,2\right\} ,+\right) $ is not $u$%
-normal, it is simply $AU$-composite.\bigskip

\noindent \textbf{Proposition 3.2.7 }\textsl{Any signature- or similar-prime
groupoid is }$u$\textsl{-normal.}\smallskip

\noindent \textit{Proof. }The proof is straightforward and we omit it.

\begin{flushright}
$\blacksquare $
\end{flushright}

\noindent \textbf{Proposition 3.2.8 }\textsl{The right-zero-semigroup on }$X$%
\textsl{\ is similar-prime.}\smallskip

\noindent \textit{Proof}. Let $\left( X,\bullet \right) $ be the
right-zero-semigroup on $X$. Then $x\bullet y=y$ for all $x,y\in X$. Let $%
\left( X,\ast \right) =U\left( X,\bullet \right) $ and $\left( X,\circ
\right) =A\left( X,\bullet \right) $, thus%
\begin{equation*}
x\ast y=%
\begin{cases}
x & \text{ if }x=y, \\ 
x\bullet y=y & \text{otherwise}%
\end{cases}%
\quad \text{and}\quad x\circ y=%
\begin{cases}
x\bullet x=x & \text{ if }x=y \\ 
x\circ y=x & \text{otherwise}%
\end{cases}%
\end{equation*}%
\noindent Hence for all $x,y$ $\in X$, $\left( X,\bullet \right) =(X,\bullet
)\diamond id_{Bin\left( X\right) }$.

\begin{flushright}
$\blacksquare $
\end{flushright}

\noindent \textbf{Example 3.2.9} Let $\left( X,\bullet \right)
=(\{a,b,c\},\bullet )$ be the right-zero-semigroup on $\left\{ a,b,c\right\} 
$. Its Cayley table together with its associated \textit{signature}-\textit{%
similar}-product tables, respectively, are:%
\begin{equation*}
\begin{tabular}{l|lll}
$\bullet $ & $a$ & $b$ & $c$ \\ \hline
$a$ & $a$ & $b$ & $c$ \\ 
$b$ & $a$ & $b$ & $c$ \\ 
$c$ & $a$ & $b$ & $c$%
\end{tabular}%
\end{equation*}%
\begin{equation*}
\text{\quad }%
\begin{tabular}{l|lll}
$\ast $ & $a$ & $b$ & $c$ \\ \hline
$a$ & $a$ & $b$ & $c$ \\ 
$b$ & $a$ & $b$ & $c$ \\ 
$c$ & $a$ & $b$ & $c$%
\end{tabular}%
\text{\quad }\diamond \text{\quad }%
\begin{tabular}{l|lll}
$\circ $ & $a$ & $b$ & $c$ \\ \hline
$a$ & $a$ & $a$ & $a$ \\ 
$b$ & $b$ & $b$ & $b$ \\ 
$c$ & $c$ & $c$ & $c$%
\end{tabular}%
\text{\quad }=\text{\quad }%
\begin{tabular}{l|lll}
$\bullet $ & $a$ & $b$ & $c$ \\ \hline
$a$ & $a$ & $b$ & $c$ \\ 
$b$ & $a$ & $b$ & $c$ \\ 
$c$ & $a$ & $b$ & $c$%
\end{tabular}%
\text{\quad }
\end{equation*}%
\noindent Therefore, the right-zero-semigroup of order 3 is \textit{similar}%
-prime since its \textit{similar}-factor $A(\{a,b,c\},\bullet )$ is $%
id_{Bin\left( X\right) }$, i.e., the left-zero-semigroup for $\left\{
a,b,c\right\} $.\medskip

\noindent \textbf{Proposition 3.2.10 }\textsl{A non-locally-zero strong
groupoid is }$u$\textsl{-composite}.\smallskip

\noindent \textit{Proof}. Let $\left( X,\bullet \right) $ $\in
Bin(X)-ZBin\left( X\right) $, then $x\bullet y\neq \{x,y\}$ for any $x,y\in
X $. Meaning, $\left( X,\bullet \right) $ cannot be the left- nor the
right-zero-semigroup on $X$. By Proposition 3.2.5, $\left( X,\bullet \right) 
$ is $u$-normal. Let $\left( X,\ast \right) =U\left( X,\bullet \right) $ and 
$\left( X,\circ \right) =A\left( X,\bullet \right) $, then%
\begin{equation*}
x\ast y=%
\begin{cases}
x & \text{ if }x=y, \\ 
x\bullet y & \text{otherwise}%
\end{cases}%
\quad \text{and}\quad x\circ y=%
\begin{cases}
x\bullet x & \text{ if }x=y \\ 
x\circ y=x & \text{otherwise}%
\end{cases}%
\end{equation*}%
\noindent Hence, for all $x,y$ $\in X$, $\left( X,\ast \right) \neq
(X,\bullet )\neq \left( X,\circ \right) $ and $\left( X,\ast \right) \neq
id_{Bin\left( X\right) }\neq \left( X,\circ \right) $. Therefore, $\left(
X,\bullet \right) $ is $u$-composite.

\subsection{Factoring $\mathbf{U}\left( X,\bullet \right) $ and $\mathbf{A}%
\left( X,\bullet \right) $}

Let $Str\left( X\right) $ be the collection of all strong groupoids on a
non-empty set $X$. Consider a groupoid $\left( X,\bullet \right) \in
Str\left( X\right) $, we classify the \textit{signature}- and \textit{similar%
}-factors of $\left( X,\bullet \right) $ as $UA$-composite, \textit{signature%
}- or \textit{similar}-prime. We conclude that $U\left( X,\bullet \right) $
and $A\left( X,\bullet \right) $ are \textit{similar}- and \textit{signature}%
-prime, respectively.\bigskip

\noindent \textbf{Theorem 3.3.1 }\textsl{The signature-factor of a strong
groupoid is similar-prime, and the similar-factor is signature-prime.}%
\smallskip

\noindent \textit{Proof.} Let $\left( X,\bullet \right) \in Str\left(
X\right) $. Suppose that $\left( X,\ast \right) =U\left( X,\bullet \right) $
and $\left( X,\circ \right) =A\left( X,\bullet \right) $. Let $\left(
X,\circledast \right) =U\left( X,\ast \right) $ and $\left( X,\odot \right)
=A\left( X,\ast \right) $, then \textquotedblleft $\circledast $%
\textquotedblright\ and \textquotedblleft $\odot $\textquotedblright\ are
defined as:%
\begin{equation*}
x\circledast y=%
\begin{cases}
x; & \text{if $x=y$,} \\ 
x\ast y=x\bullet y & \text{otherwise}%
\end{cases}%
\text{\quad and\quad }x\odot y=%
\begin{cases}
x\ast x=x & \text{if $x=y$,} \\ 
x; & \text{otherwise.}%
\end{cases}%
\end{equation*}%
\smallskip

\noindent Hence $A\left( X,\ast \right) =id_{Bin\left( X\right) }$, and
therefore $U\left( X,\bullet \right) $ is \textit{similar}-prime. Similarly,
if we let $\left( X,\boxtimes \right) =U\left( X,\circ \right) $ and $\left(
X,\boxdot \right) =A\left( X,\circ \right) $, then \textquotedblleft $%
\boxtimes $\textquotedblright\ and \textquotedblleft $\boxdot $%
\textquotedblright\ are defined as:

\begin{equation*}
x\boxtimes y=%
\begin{cases}
x & \text{if $x=y$,} \\ 
x\circ y=x & \text{otherwise}%
\end{cases}%
\text{\quad and\quad }\ x\boxdot y=%
\begin{cases}
x\circ x=x\bullet x & \text{if $x=y$,} \\ 
x; & \text{otherwise.}%
\end{cases}%
\end{equation*}%
\smallskip

\noindent Therefore, $U\left( X,\circ \right) =id_{Bin\left( X\right) }$,
and hence $A\left( X,\bullet \right) $ is \textit{signature}-prime.

\begin{flushright}
$\blacksquare $
\end{flushright}

\noindent \textbf{Corollary 3.3.2.} \textbf{\ }\textsl{Let }$\left(
X,\bullet \right) $\textsl{\ be any groupoid and let }$\left( X,\ast \right)
=U\left( X,\bullet \right) $\textsl{\ and }$\left( X,\circ \right) =A\left(
X,\bullet \right) $\textsl{. If }$\left( X,\bullet \right) $\textsl{\ has a }%
$UA$\textsl{-factorization, i.e., if }$\left( X,\bullet \right) =\left(
X,\ast \right) \diamond \left( X,\circ \right) ,$\textsl{\ then }%
\begin{equation*}
\left( X,\bullet \right) =U\left( X,\ast \right) \diamond A\left( X,\circ
\right) .
\end{equation*}

\noindent \textit{Proof.} This follows immediately from the previous
theorem. In fact, suppose $\left( X,\bullet \right) $ has a $UA$%
-factorization, then%
\begin{eqnarray*}
\left( X,\bullet \right) &=&\left( X,\ast \right) \diamond \left( X,\circ
\right) \\
&=&(U\left( X,\ast \right) \diamond A\left( X,\ast \right) )\diamond \left(
U\left( X,\circ \right) \diamond A\left( X,\circ \right) \right) \\
&=&(U\left( X,\ast \right) \diamond id_{Bin\left( X\right) })\diamond \left(
id_{Bin\left( X\right) }\diamond A\left( X,\circ \right) \right) \\
&=&U\left( X,\ast \right) \diamond A\left( X,\circ \right) \text{.}
\end{eqnarray*}

\begin{flushright}
$\blacksquare $
\end{flushright}

\noindent \textbf{Corollary 3.3.3.} \textbf{\ }\textsl{Let }$\left(
X,\bullet \right) $\textsl{\ be a groupoid and let }$\left( X,\ast \right)
=U\left( X,\bullet \right) $\textsl{\ and }$\left( X,\circ \right) =A\left(
X,\bullet \right) $\textsl{. If }$\left( X,\bullet \right) $\textsl{\ has a }%
$AU$\textsl{-factorization then }%
\begin{equation*}
\left( X,\bullet \right) =A\left( X,\circ \right) \diamond U\left( X,\ast
\right) .
\end{equation*}

\noindent \textit{Proof.} The proof is very similar to that of the previous
Corollary.

\begin{flushright}
$\blacksquare $
\end{flushright}

\noindent \textbf{Corollary 3.3.4.\ }\textsl{Let }$\left( X,\bullet \right) $%
\textsl{\ be a strong groupoid and let }$\left( X,\ast \right) =U\left(
X,\bullet \right) $\textsl{\ and }$\left( X,\circ \right) =A\left( X,\bullet
\right) $\textsl{, then }%
\begin{equation*}
\left( X,\bullet \right) =A\left( X,\circ \right) \diamond U\left( X,\ast
\right) =U\left( X,\ast \right) \diamond A\left( X,\circ \right) \text{.}
\end{equation*}%
\smallskip \textsl{\ }\noindent \textit{Proof.} This is a direct result of
Theorem 3.1.3 and the previous two Corollaries.

\begin{flushright}
$\blacksquare $
\end{flushright}

As a final observation, a groupoid is \textit{similar}-prime if it is
similar to the left-zero-semigroup or a locally-zero-groupoid, in other
words, if it is idempotent. Hence, we need another method of factorization
for idempotent groupoids.

\section{\protect\bigskip Orient-Skew Factorization}

We say a groupoid $\left( X,\ast \right) \,$\ has the \emph{orientation
property OP }\cite{KN08} if $x\ast y\in \{x,y\}$ for all $x,y\in X$.
Moreover, $\left( X,\ast \right) $ has the \emph{twisted orientation
property TOP} if $x\ast y=x$ implies $y\ast x=x$ for all $x,y\in X$. In this
section, we introduce a unique factorization which can be applied to
groupoids with \emph{OP}. This type of groupoids has proven to be useful in
graph theory, where in a directed graph $x\ast y=x$ can mean there is a path
from vertex $x$ to vertex $y$, i.e. $x\rightarrow y$; while $x\ast y=y$ can
mean there is no path from $x$ to $y$, i.e. $x\nrightarrow y$. In fact, if $%
\Gamma _{\left( X,\ast \right) }$ is the directed graph on vertex set $X$
and $\left( X,\ast \right) \in TOP\left( X\right) $, then $\Gamma _{\left(
X,\ast \right) }$ is a simple graph \cite{AKN18}. For more details on
groupoids associated with directed and simple graphs we refer to \cite%
{AKN18, We01}.\bigskip

\noindent \textbf{Example 4.1 }Let $X=\{0,1\}$ and $\left( X,\leq \right) $
be a linearly ordered set. Define a binary operation \textquotedblleft $%
\bullet $\textquotedblright\ on $X$ such that:%
\begin{equation*}
x\bullet y=%
\begin{cases}
0 & \text{ if $x\leq y,$} \\ 
1 & \text{otherwise.}%
\end{cases}%
\end{equation*}%
\noindent Then the binary system $\left( X,\bullet \right) $ has the
orientation property.\medskip

\noindent \textbf{Example 4.2} Let $X=\left\{ a,b,c\right\} $. Define a
binary operation \textquotedblleft $\bullet $\textquotedblright\ on $X$ by
the following table:%
\begin{equation*}
\begin{tabular}{r|rrr}
$\bullet $ & $a$ & $b$ & $c$ \\ \hline
$a$ & $a$ & $b$ & $c$ \\ 
$b$ & $b$ & $b$ & $c$ \\ 
$c$ & $c$ & $b$ & $c$%
\end{tabular}%
\end{equation*}%
\noindent Then $\left( X,\bullet \right) $ has the twisted orientation
property.\bigskip

We consider three functions to represent operations on the main diagonal and
on the anti-diagonal of the associated Cayley table of a binary operation on
a finite set.

Let $(X,\ast )$ be a groupoid of finite order $n$ and binary operation
\textquotedblleft $\ast $\textquotedblright , i.e., $|X|=n$ and $\ast
:X^{2}\longrightarrow X$. Then for all $x_{i}$, $x_{j}\in X$, $i,j=1,2,...,n$%
, and $i+j=n+1$, we call:

\begin{description}
\item[diag-1] $\overline{d^{\ast }}$ the \emph{anti-diagonal function} of $%
(X,\ast )\ $such that $\overline{d^{\ast }}$ $:%
\mathbb{N}
\longrightarrow X$, defined by $\overline{d^{\ast }}$ $\left( i\right)
=x_{i}\ast x_{j}.$

\item[diag-2] $\widehat{d^{\ast }}$ the \emph{reverse-diagonal function} of $%
(X,\ast )\ $such that $\widehat{d^{\ast }}$ $:%
\mathbb{N}
\longrightarrow X$, defined by $\widehat{d^{\ast }}$ $\left( i\right)
=x_{j}\ast x_{j}.$

\item[diag-3] $\widetilde{d^{\ast }}$ the \emph{skew-diagonal function} of $%
(X,\ast )\ $such that $\widetilde{d^{\ast }}$ $:%
\mathbb{N}
\longrightarrow X$, defined by $\widetilde{d^{\ast }}\left( i\right) =%
\widehat{\overline{d^{\ast }}}(i)=$ $x_{j}\ast x_{i}.$\bigskip
\end{description}

\noindent \textbf{Example 4.3 }Consider the groupoid $\left( \left\{
0,1,2,3\right\} ,\ast \right) $ where \textquotedblleft $\ast $%
\textquotedblright\ is given by the following table:%
\begin{equation*}
\begin{tabular}{c|cccc}
$\ast $ & 0 & 1 & 2 & 3 \\ \hline
0 & 0 & 1 & 0 & 3 \\ 
1 & 1 & 1 & 1 & 0 \\ 
2 & 2 & 2 & 2 & 3 \\ 
3 & 0 & 3 & 2 & 3%
\end{tabular}%
\end{equation*}%
\noindent Observe that $n=4$ and the main diagonal $d^{\ast }=\{0,1,2,3\}$.
For instance, $d^{\ast }(2)=2\ast 2=2$. Also, the anti-diagonal $\overline{%
d^{\ast }}=\{3,1,2,0\}$. For example, $\overline{d^{\ast }}\left( 1\right)
=x_{1}\ast x_{4}=0\ast 3=3$. Moreover, the reverse of the diagonal is $%
\widehat{d^{\ast }}=\left\{ 3,2,1,0\right\} $. For instance, $\widehat{%
d^{\ast }}$ $\left( 4\right) =x_{1}\ast x_{1}=0\ast 0=0$. So the
skew-diagonal defined here is the reverse of the anti-diagonal, hence, $%
\widetilde{d^{\ast }}$ $=\left\{ 0,2,1,3\right\} $. For example, $\widetilde{%
d^{\ast }}\left( 3\right) =\widehat{\overline{d^{\ast }}}(3)=$ $x_{2}\ast
x_{3}=1\ast 2=1$.\medskip

Given these definitions, we can \textit{derive} the \textit{orient}-factor
of a groupoid from $id_{Bin\left( X\right) }$, such that all its elements
are the same as those of the left-zero-semigroup except elements belonging
to the anti-diagonal, which we construct from the skew-diagonal of $%
id_{Bin\left( X\right) }$. Similarly, the \textit{skew}-factor is \textit{%
derived} from the parent groupoid by letting its anti-diagonal be that of
the skew-diagonal of the parent groupoid, otherwise all other elements are
kept the same as the parent groupoid. \medskip

Let $\left( X,\bullet \right) $ be a groupoid. Let $D^{\diamond }$ denote
the main diagonal of $id_{Bin\left( X\right) }$. \textit{Derive} groupoids $%
(X,\ast )$ and $\left( X,\circ \right) $ from $id_{Bin\left( X\right) }$ and 
$\left( X,\bullet \right) $, respectively, as follows:

\noindent For all $x,y\in X$,%
\begin{equation}
\begin{tabular}{lll}
$\text{(i) }\overline{d^{\ast }}=\widetilde{D^{\diamond }}$, & and$\text{%
\quad }$ & (i) $\overline{d^{\circ }}=\widetilde{d^{\bullet }}$, \\ 
$\text{(ii) }x\ast y=x\text{; otherwise.}$ &  & (ii) $x\circ y=x\bullet y$;
otherwise.%
\end{tabular}
\label{Orient-Skew-F's}
\end{equation}%
Groupoids $\left( X,\ast \right) $ and $\left( X,\circ \right) $ are said to
be the \emph{orient-} and \emph{skew-factor }of $\left( X,\bullet \right) $,
respectively, denoted by $O\left( X,\bullet \right) $ and $J\left( X,\bullet
\right) $. As previously mentioned, the product \textquotedblleft $\diamond $%
\textquotedblright\ is not commutative. Hence, for $\left( X,\bullet \right)
\in Bin\left( X\right) $, we may have an $OJ$\emph{-factorization} such that%
\begin{equation}
\left( X,\bullet \right) =O\left( X,\bullet \right) \diamond J\left(
X,\bullet \right)  \label{OJ-prod}
\end{equation}%
or a $JO$\emph{-factorization} such that%
\begin{equation}
\left( X,\bullet \right) =J\left( X,\bullet \right) \diamond O\left(
X,\bullet \right) \text{.}  \label{JO-prod}
\end{equation}%
\bigskip

\noindent \textbf{Proposition 4.4 }\textsl{The orient-factor of a given
groupoid is locally-zero.\smallskip }

\noindent \textit{Proof}. Given $\left( X,\bullet \right) \in Bin\left(
X\right) $, let $(X,\ast )=O\left( X,\bullet \right) $. Then, $d^{\ast
}=D^{\diamond }$, i.e. $x\ast x=x$, and $x\ast y=x$ for all $x,$ $y\in X$
except when $x,$ $y\in \overline{d^{\ast }}$. In fact, for any $x\not=y$ in $%
X$, $(\{x,y\},\bullet )$ is either a left-zero-semigroup or a
right-zero-semigroup. Moreover, $x\bullet x=x$ for all $x\in X$ which
implies that $O\left( X,\bullet \right) $ is locally-zero.

\begin{flushright}
$\blacksquare $
\end{flushright}

\noindent \textbf{Corollary 4.5 }\textsl{The orient-factor of a given
groupoid is a unit in }$Bin\left( X\right) $\textsl{.\smallskip }

\noindent \textit{Proof}. This follows immediately from Propositions 2.8 and
4.4.

\begin{flushright}
$\blacksquare $
\end{flushright}

\noindent \textbf{Example 4.6 }Let $X=\{e,a,b,c\}$. Define a binary
operation \textquotedblleft $\bullet $\textquotedblright\ by the following
table:%
\begin{equation*}
\begin{tabular}{c|cccc}
$\bullet $ & $e$ & $a$ & $b$ & $c$ \\ \hline
$e$ & $e$ & $a$ & $b$ & $c$ \\ 
$a$ & $a$ & $e$ & $c$ & $b$ \\ 
$b$ & $b$ & $c$ & $a$ & $e$ \\ 
$c$ & $c$ & $b$ & $e$ & $a$%
\end{tabular}%
\ 
\end{equation*}%
\noindent Then, clearly $\left( X,\bullet ,e\right) $ is a group. Derive its 
\textit{orient-}factor $\left( X,\ast ,e\right) =U(X,\bullet ,e)$ as in \ref%
{Orient-Skew-F's} to obtain:%
\begin{equation*}
\begin{tabular}{c|cccc}
$\ast $ & $e$ & $a$ & $b$ & $c$ \\ \hline
$e$ & $e$ & $e$ & $e$ & $c$ \\ 
$a$ & $a$ & $a$ & $b$ & $a$ \\ 
$b$ & $b$ & $a$ & $b$ & $b$ \\ 
$c$ & $e$ & $c$ & $c$ & $c$%
\end{tabular}%
\end{equation*}%
\noindent Hence, $\left( X,\ast ,e\right) $ is locally-zero.

\subsection{$OJ$-Factorization}

In this subsection, we explore an $OJ$-factorization of any groupoid $\left(
X,\bullet \right) $ in $Bin\left( X\right) $, i.e., into its \textit{orient-}
and \textit{skew}-factors, respectively. The next subsection discusses a $JO$%
-factorization where the product of the two factors is \textquotedblleft
reversed\textquotedblright . Then, we classify $\left( X,\bullet \right) $
as $OJ$\textit{- }and/or\textit{\ }$JO$-\textit{composite}, $j$-\textit{%
composite} or\textit{\ }$j$\textit{-normal;\textit{\ }}and as\textit{\
orient- }or\textit{\ skew-prime.}\medskip

A groupoid $\left( X,\bullet \right) $ is \emph{bi-diagonal} if its
anti-diagonal is symmetric, meaning if $\overline{d^{\bullet }}=\widetilde{%
d^{\bullet }}$.\medskip

\noindent \textbf{Example 4.1.1.} Let $\,(%
\mathbb{Z}
,<)$ be a linearly ordered set. Consider groupoid $\left( 
\mathbb{Z}
,\bullet \right) $ where $x\bullet y=$ $\max \left\{ x,y\right\} $ for all $%
x,y\in 
\mathbb{Z}
$. Define two binary operations on $%
\mathbb{Z}
$ such that: 
\begin{equation*}
x\ast y=%
\begin{cases}
x & \text{if }x<y, \\ 
y & \text{otherwise.}%
\end{cases}%
\text{\quad and\quad }x\circ y=%
\begin{cases}
x & \text{if $x\leq y$,} \\ 
y & \text{otherwise.}%
\end{cases}%
\noindent
\end{equation*}%
\noindent Then clearly $(X,\ast )\diamond (X,\circ )$\ is an $OJ$%
-factorization of $\left( X,\bullet \right) $, where $(X,\ast )=O(X,\bullet
) $ and $(X,\circ )=J(X,\bullet )$. Moreover, $(%
\mathbb{Z}
,$ $\bullet )$ is bi-diagonal.\medskip

A groupoid $\left( X,\bullet \right) $ is said to be \emph{orient-prime}%
\textit{\ }if $O\left( X,\bullet \right) $ $=id_{Bin\left( X\right) }$, and
is said to be \emph{skew-prime} if $J\left( X,\bullet \right) $ $%
=id_{Bin\left( X\right) }$. Alternatively, if $\left( X,\bullet \right) $ is
neither \textit{orient}- nor \textit{skew}-prime, then $\left( X,\bullet
\right) $ is said to be

\qquad (1) $OJ$\emph{-composite} if $\left( X,\bullet \right) =O\left(
X,\bullet \right) \diamond J\left( X,\bullet \right) $;

\qquad (2) $JO$\emph{-composite} if $\left( X,\bullet \right) =J\left(
X,\bullet \right) \diamond O\left( X,\bullet \right) $.

Consequently, $\left( X,\bullet \right) $ is said to be $j$\emph{-composite}
if both $\left( 1\right) $ and $\left( 2\right) $ hold.\medskip

Just as with $UA$-factorization, not every groupoid will have a $JO$%
-factorization. But it is possible to \textit{derive} an $OJ$-factorization
of any given groupoid.\medskip\ 

\noindent \textbf{Theorem 4.1.2} \textsl{Any given groupoid has an }$OJ$%
\textsl{-factorization, i.e., if }$(X,\bullet )\in Bin(X)$\textsl{$,$ then}%
\begin{equation*}
\left( X,\bullet \right) =O(X,\bullet )\diamond J(X,\bullet )\text{.}
\end{equation*}%
\noindent \textit{Proof.} Let $(X,\bullet )$ $\in Bin\left( X\right) $ such
that $O(X,\bullet )$ and $J(X,\bullet )$ are defined as in \ref%
{Orient-Skew-F's}. Let $(X,\odot )=(X,\ast )$ $\diamond (X,\circ )$ where $%
(X,\ast )=O(X,\bullet )$ and $(X,\circ )=J(X,\bullet )$. Then $x\odot
y=(x\ast y)\circ (y\ast x)$ for all $x,y\in X$. It follows that

\begin{enumerate}
\item[(i)] If $x=y,$ $x\ast x=x$ and $x\circ x=x\bullet x$.

\item[(ii)] If $x\not=y$, then if $x\ast y\in \overline{d^{\ast }},$ $x\ast
y\in $ $\widetilde{D^{\diamond }}$, and for $x\circ y\in \overline{d^{\circ }%
},$ then $x\circ y\in $ $\widetilde{d^{\bullet }}$. Otherwise, $x\ast y=x$,
and $x\circ y=x\bullet y$.
\end{enumerate}

Next, we show that $(X,\bullet )=(X,\odot )$. Given $x,y\in X$,

\begin{enumerate}
\item[(i)] If $x=y,$ $x\odot x=(x\circ x)\ast (x\circ x)=x\circ x=x\bullet x$%
.

\item[(ii)] If $x\not=y$, then if $x\ast y=y\ast x,$ then $x\odot y=(x\ast
y)\circ (y\ast x)=x\circ y=x\bullet y$ and $y\odot x=(y\ast x)\circ (x\ast
y)=y\circ x=x\bullet y$. If $x\ast y\not=y\ast x,$ then $x\odot y=(x\ast
y)\circ (y\ast x)=(x\ast y)\bullet (y\ast x)\in \left\{ x\bullet y,\text{ }%
y\bullet x\right\} $.
\end{enumerate}

\noindent Thus, $x\odot y=x\bullet y$ for all $x,y\in X$. This proves that $%
(X,\odot )=(X,\bullet )$.

\begin{flushright}
$\blacksquare $
\end{flushright}

\noindent \textbf{Corollary 4.1.3} \textsl{The factorization in Theorem
4.1.2 is unique.} \smallskip

\noindent \textit{Proof.} Let $(X,\bullet )$ $\in Bin\left( X\right) $ with
an $OJ$-factorization such that $(X,\bullet )=(X,\ast )$ $\diamond (X,\circ
) $ where $(X,\ast )=O(X,\bullet )$ and $(X,\circ )=J(X,\bullet )$. Let $%
(X,\bullet )=(X,\bigtriangledown )$ $\diamond (X,\bigtriangleup )$ where $%
(X,\bigtriangledown )=O(X,\bullet )$ and $(X,\bigtriangleup )=J(X,\bullet )$%
. For any $x\in X$, we have $x\ast x=x=x\bigtriangledown x$, and $x\ast
y=x\bigtriangledown y$ when $x\not=y$. Hence $(X,\ast )=(X,\bigtriangledown
) $. Similarly, if $x\in X$, then $x\circ x=x\bullet x=x\bigtriangleup x$.
When $x\not=y$, we have $x\circ y=x\bullet y=x\bigtriangleup y$, proving
that $(X,\circ )=(X,\bigtriangleup )$.

\begin{flushright}
$\blacksquare $
\end{flushright}

\noindent \textbf{Example 4.1.4} \cite{NJK99} Consider the groupoid $\left(
X,\bullet \right) =\left( \left\{ 1,2,3,4\right\} ,\bullet \right) $ where
\textquotedblleft $\bullet $\textquotedblright\ is defined by the following
Cayley table:%
\begin{equation*}
\begin{tabular}{c|cccc}
$\bullet $ & 1 & 2 & 3 & 4 \\ \hline
1 & 1 & 1 & 3 & 1 \\ 
2 & 2 & 2 & 3 & 2 \\ 
3 & 1 & 2 & 3 & 4 \\ 
4 & 4 & 4 & 3 & 4%
\end{tabular}%
\ 
\end{equation*}%
\noindent By \textit{deriving} its \textit{orient}- and \textit{skew}%
-factors $O\left( X,\bullet \right) $ and $J\left( X,\bullet \right) $,
respectively, and by letting $\left( X,\ast \right) =O\left( X,\bullet
\right) $ and $\left( X,\circ \right) =J\left( X,\bullet \right) $ shows
that $\left( X,\ast \right) \diamond \left( X,\circ \right) =$\ $\left(
X,\bullet \right) $.

\noindent Indeed, $\left( X,\bullet \right) $ has an $OJ$-factorization:%
\textsl{\smallskip }%
\begin{equation*}
\text{\quad }%
\begin{tabular}{c|cccc}
$\ast $ & 1 & 2 & 3 & 4 \\ \hline
1 & 1 & 1 & 1 & 4 \\ 
2 & 2 & 2 & 3 & 2 \\ 
3 & 3 & 2 & 3 & 3 \\ 
4 & 1 & 4 & 4 & 4%
\end{tabular}%
\text{\quad }\diamond \text{\quad }%
\begin{tabular}{c|cccc}
$\circ $ & 1 & 2 & 3 & 4 \\ \hline
1 & 1 & 1 & 3 & 4 \\ 
2 & 2 & 2 & 2 & 2 \\ 
3 & 1 & 3 & 3 & 4 \\ 
4 & 1 & 4 & 3 & 4%
\end{tabular}%
\text{\quad }=\text{\quad }%
\begin{tabular}{c|cccc}
$\bullet $ & 1 & 2 & 3 & 4 \\ \hline
1 & 1 & 1 & 3 & 1 \\ 
2 & 2 & 2 & 3 & 2 \\ 
3 & 1 & 2 & 3 & 4 \\ 
4 & 4 & 4 & 3 & 4%
\end{tabular}%
\end{equation*}%
\noindent Also, since $O\left( X,\bullet \right) \neq id_{Bin\left( X\right)
}\neq J\left( X,\bullet \right) $, then $\left( X,\bullet \right) $ is $OJ$%
-composite.\medskip

\subsection{$JO$-Factorization}

In this subsection, we reverse the product of the \textit{orient}- and 
\textit{skew}-factors of a given groupoid $\left( X,\bullet \right) \in
Bin\left( X\right) $. We find that an arbitrary groupoid admits a $JO$%
-factorization if it has the orientation property.\medskip

\noindent \textbf{Example 4.2.1 }Consider the groupoid $\left( X,\bullet
\right) =\left( \left\{ 1,2,3,4\right\} ,\bullet \right) $ defined as in
Example 4.1.4:%
\begin{equation*}
\begin{tabular}{c|cccc}
$\bullet $ & 1 & 2 & 3 & 4 \\ \hline
1 & 1 & 1 & 3 & 1 \\ 
2 & 2 & 2 & 3 & 2 \\ 
3 & 1 & 2 & 3 & 4 \\ 
4 & 4 & 4 & 3 & 4%
\end{tabular}%
\ 
\end{equation*}%
\noindent Through routine calculations, we find that $\left( X,\bullet
\right) $ admits a $JO$-factorization since $J\left( X,\bullet \right)
\diamond O\left( X,\bullet \right) =\left( X,\circ \right) \diamond \left(
X,\ast \right) =$\ $\left( X,\bullet \right) .$ In addition, $\left(
X,\bullet \right) \in OP\left( X\right) $.\medskip

A groupoid $\left( X,\bullet \right) $ is said to be $j$\emph{-normal}%
\textit{\ }if it admits an $OJ$- and a $JO$-factorization, i.e., if

\qquad (i) $\left( X,\bullet \right) =O\left( X,\bullet \right) \diamond
J\left( X,\bullet \right) $ and

\qquad (ii) $\left( X,\bullet \right) =J\left( X,\bullet \right) \diamond
O\left( X,\bullet \right) $.\bigskip

\noindent \textbf{Theorem 4.2.3} \textsl{A groupoid }$(X,\bullet )$\textsl{\
with the orientation property has a }$JO$\textsl{-factorization.\smallskip }

\noindent \textit{Proof.} Let $(X,\bullet )$ $\in OP\left( X\right) $.
Define $(X,\odot )=(X,\ast )$ $\diamond (X,\circ )$ where $(X,\ast
)=O(X,\bullet )$ and $(X,\circ )=J(X,\bullet )$. Then $x\odot y=(x\ast
y)\circ (y\ast x)$ for all $x,y\in X$. It follows that

\begin{enumerate}
\item[(i)] If $x=y,$ then $x\ast x=x$ and $x\circ x=x\bullet x.$

\item[(ii)] If $x\not=y,$ the two cases arise: if $x\ast y\in \overline{%
d^{\ast }}$ and $x\circ y\in \overline{d^{\circ }},$ then $x\ast y\in $ $%
\widetilde{D^{\diamond }}$ and $x\circ y\in $ $\widetilde{d^{\bullet }}$
which also $\in \left\{ x,\text{ }y\right\} $. Otherwise, $x\ast y=x,$and $%
x\circ y=x\bullet y.$
\end{enumerate}

\noindent Next, we show that $(X,\bullet )=(X,\odot )$. Given $x,y\in X$,

\begin{enumerate}
\item[(i)] If $x=y,$ then $x\odot x=(x\circ x)\ast (x\circ x)=x\ast
x=x\bullet x$.

\item[(ii)] If $x\not=y,$ then $x\odot y=(x\circ y)\ast (y\circ x)$. If $%
x\circ y=y\circ x$, then $x\odot y=(x\circ y)\ast (x\circ y)=x\circ
y=x\bullet y$. If $x\circ y\not=y\circ x$, then $x\odot y=(x\ast y)\bullet
(y\ast x)\in \left\{ x\bullet y,\text{ }y\bullet x\right\} $.
\end{enumerate}

\noindent Thus $x\odot y=x\bullet y$ for all $x,y\in X$. This proves that $%
(X,\odot )=(X,\bullet )$.

\begin{flushright}
$\blacksquare $
\end{flushright}

\noindent \textbf{Corollary 4.2.4} \textsl{The factorization in Theorem
4.2.3 is unique.} \smallskip

\noindent \textit{Proof. }The proof is very similar to that of Corollary
4.1.3 so we omit it.

\begin{flushright}
$\blacksquare $
\end{flushright}

\noindent \textbf{Proposition 4.2.5} \textsl{A groupoid with the orientation
property is }$j$\textsl{-normal.\smallskip }

\noindent \textit{Proof.} The result follows from Theorems 4.1.2, 4.2.3 and
the definition.

\begin{flushright}
$\blacksquare $
\end{flushright}

\noindent \textbf{Example 4.2.6} Let $\left( X,\bullet \right) $ be defined
as in Example 4.2.1 where we determined that $\left( X,\bullet \right) $
admits an $OJ$-factorization. It can be verified that $J\left( X,\bullet
\right) \diamond O\left( X,\bullet \right) =\left( X,\bullet \right) $,
which shows that $\left( X,\bullet \right) $ admits a $JO$-factorization as
well. Therefore, $\left( X,\bullet \right) $ is $j$-normal in $(Bin\left(
X\right) ,\diamond )$. Additionally, $J\left( X,\bullet \right) \neq
id_{Bin\left( X\right) }\neq O\left( X,\bullet \right) $ implies that $%
\left( X,\bullet \right) $ is $j$-composite.

\subsection{Factoring $\mathbf{O}\left( X,\bullet \right) $ and $\mathbf{J}%
\left( X,\bullet \right) $}

In this subsection, the \textit{orient}- and \textit{skew}-factors of $%
\left( X,\bullet \right) \in OP\left( X\right) $ are factored to deduce that 
$O\left( X,\bullet \right) $ is \textit{skew}-prime while $J\left( X,\bullet
\right) $ is \textit{binary}-equivalent to $\left( X,\bullet \right) $%
.\medskip

Let $\left( X,\bullet \right) $ and $\left( X,\circ \right) $ be groupoids
in $Bin\left( X\right) $. We say that $\left( X,\circ \right) $ is \emph{%
binary-equivalent} to $\left( X,\bullet \right) $ if there exists $\left(
X,\ast \right) \in Bin\left( X\right) $ such that

$\qquad $(i) $\left( X,\bullet \right) =\left( X,\ast \right) \diamond
\left( X,\circ \right) $; and

$\qquad $(ii) $\left( X,\circ \right) =\left( X,\ast \right) \diamond \left(
X,\bullet \right) $.\medskip

\noindent \textbf{Theorem 4.3.1 }\textsl{Given a groupoid }$\left( X,\bullet
\right) $\textsl{\ with the orientation property. Its orient-factor\ is
skew-prime, and its skew-factor is binary-equivalent to }$\left( X,\bullet
\right) $\textsl{.\smallskip }

\noindent \textit{Proof.} Let $\left( X,\bullet \right) \in OP\left(
X\right) $. Suppose that $\left( X,\ast \right) =O\left( X,\bullet \right) $
and $\left( X,\circ \right) =J\left( X,\bullet \right) $. Then by Theorem
4.1.2 $\left( X,\bullet \right) =O(X,\bullet )\diamond J(X,\bullet )=\left(
X,\ast \right) \diamond \left( X,\circ \right) $. Let $\left( X,\circledast
\right) =O\left( X,\ast \right) $ and $\left( X,\odot \right) =J\left(
X,\ast \right) ,$ then for $\circledast $: (i)$\overline{\text{ }%
d^{\circledast }}=\widetilde{D^{\diamond }}$, (ii) $x\circledast y=x$,
otherwise; and for $\odot $: (i) $\overline{d^{\odot }}=\widetilde{d^{\ast }}%
=D^{\diamond }$, (ii) $x\odot y=x\ast y=x$, otherwise. Hence,%
\begin{equation*}
\left( X,\ast \right) =\left( X,\ast \right) \diamond id_{Bin\left( X\right)
}
\end{equation*}

\noindent and $O\left( X,\bullet \right) $ is \textit{skew}-prime.
Similarly, if we let $\left( X,\boxtimes \right) =O\left( X,\circ \right) $
and $\left( X,\boxdot \right) =J\left( X,\circ \right) $, then for $%
\boxtimes $: (i)$\overline{\text{ }d^{\boxtimes }}=\widetilde{D^{\diamond }} 
$, (ii) $x\boxtimes y=x$, otherwise; and for $\boxdot $: (i) $\overline{%
d^{\boxdot }}=\widetilde{d^{\circ }}=\overline{d^{\bullet }}$, (ii) $%
x\boxdot y=x\circ y=x\bullet y$, otherwise. Thus,%
\begin{equation*}
\left( X,\circ \right) =\left( X,\ast \right) \diamond \left( X,\bullet
\right)
\end{equation*}

\noindent and the final result follows.

\begin{flushright}
$\blacksquare $
\end{flushright}

\noindent \textbf{Example 4.3.2} Consider the locally-zero groupoid $\left(
X,\bullet \right) =\left( \left\{ 0,1,2,3,4,5\right\} ,\bullet \right) $
where \textquotedblleft $\bullet $\textquotedblright\ is defined by the
following Cayley table:%
\begin{equation*}
\begin{tabular}{c|cccccc}
$\bullet $ & 0 & 1 & 2 & 3 & 4 & 5 \\ \hline
0 & 0 & 1 & 0 & 0 & 4 & 0 \\ 
1 & 0 & 1 & 2 & 3 & 1 & 5 \\ 
2 & 2 & 1 & 2 & 3 & 4 & 2 \\ 
3 & 3 & 1 & 2 & 3 & 3 & 3 \\ 
4 & 0 & 4 & 2 & 4 & 4 & 4 \\ 
5 & 5 & 1 & 5 & 5 & 5 & 5%
\end{tabular}%
\ 
\end{equation*}%
\noindent Since $\left( X,\bullet \right) $ has the orientation property,
then $\left( X,\bullet \right) $ is $j$-normal by Proposition 4.2.5.

\noindent Factoring its \textit{orient-} and \textit{skew-}factors $\left(
X,\ast \right) =O\left( X,\bullet \right) $ and $\left( X,\circ \right)
=J\left( X,\bullet \right) $ into their respective \textit{orient}- and 
\textit{skew}-factors, $O\left( X,\ast \right) $, $J\left( X,\ast \right) $
and $O\left( X,\circ \right) $, $J\left( X,\circ \right) $, is observed
through their respective product tables:%
\begin{eqnarray*}
\text{\quad }%
\begin{tabular}{c|cccccc}
$\ast $ & 0 & 1 & 2 & 3 & 4 & 5 \\ \hline
0 & 0 & 0 & 0 & 0 & 0 & 5 \\ 
1 & 1 & 1 & 1 & 1 & 4 & 1 \\ 
2 & 2 & 2 & 2 & 3 & 2 & 2 \\ 
3 & 3 & 3 & 2 & 3 & 3 & 3 \\ 
4 & 4 & 1 & 4 & 4 & 4 & 4 \\ 
5 & 0 & 5 & 5 & 5 & 5 & 5%
\end{tabular}%
\text{\quad } &=&\text{\quad }%
\begin{tabular}{c|cccccc}
$\ast $ & 0 & 1 & 2 & 3 & 4 & 5 \\ \hline
0 & 0 & 0 & 0 & 0 & 0 & 5 \\ 
1 & 1 & 1 & 1 & 1 & 4 & 1 \\ 
2 & 2 & 2 & 2 & 3 & 2 & 2 \\ 
3 & 3 & 3 & 2 & 3 & 3 & 3 \\ 
4 & 4 & 1 & 4 & 4 & 4 & 4 \\ 
5 & 0 & 5 & 5 & 5 & 5 & 5%
\end{tabular}%
\text{\quad }\diamond \text{\quad }%
\begin{tabular}{c|cccccc}
$\odot $ & 0 & 1 & 2 & 3 & 4 & 5 \\ \hline
0 & 0 & 0 & 0 & 0 & 0 & 0 \\ 
1 & 1 & 1 & 1 & 1 & 1 & 1 \\ 
2 & 2 & 2 & 2 & 2 & 2 & 2 \\ 
3 & 3 & 3 & 3 & 3 & 3 & 3 \\ 
4 & 4 & 4 & 4 & 4 & 4 & 4 \\ 
5 & 5 & 5 & 5 & 5 & 5 & 5%
\end{tabular}%
\text{\quad }\medskip \\
\text{\quad }%
\begin{tabular}{c|cccccc}
$\circ $ & 0 & 1 & 2 & 3 & 4 & 5 \\ \hline
0 & 0 & 1 & 0 & 0 & 4 & 5 \\ 
1 & 0 & 1 & 2 & 3 & 4 & 5 \\ 
2 & 2 & 1 & 2 & 2 & 4 & 2 \\ 
3 & 3 & 1 & 3 & 3 & 3 & 3 \\ 
4 & 0 & 1 & 2 & 4 & 4 & 4 \\ 
5 & 0 & 1 & 5 & 5 & 5 & 5%
\end{tabular}%
\text{\quad } &=&\text{\quad }%
\begin{tabular}{c|cccccc}
$\ast $ & 0 & 1 & 2 & 3 & 4 & 5 \\ \hline
0 & 0 & 0 & 0 & 0 & 0 & 5 \\ 
1 & 1 & 1 & 1 & 1 & 4 & 1 \\ 
2 & 2 & 2 & 2 & 3 & 2 & 2 \\ 
3 & 3 & 3 & 2 & 3 & 3 & 3 \\ 
4 & 4 & 1 & 4 & 4 & 4 & 4 \\ 
5 & 0 & 5 & 5 & 5 & 5 & 5%
\end{tabular}%
\text{\quad }\diamond \text{\quad }%
\begin{tabular}{c|cccccc}
$\bullet $ & 0 & 1 & 2 & 3 & 4 & 5 \\ \hline
0 & 0 & 1 & 0 & 0 & 4 & 0 \\ 
1 & 0 & 1 & 2 & 3 & 1 & 5 \\ 
2 & 2 & 1 & 2 & 3 & 4 & 2 \\ 
3 & 3 & 1 & 2 & 3 & 3 & 3 \\ 
4 & 0 & 4 & 2 & 4 & 4 & 4 \\ 
5 & 5 & 1 & 5 & 5 & 5 & 5%
\end{tabular}%
\text{\quad }
\end{eqnarray*}%
\noindent Indeed, $\left( X,\ast \right) =O\left( X,\ast \right) \diamond
J\left( X,\ast \right) =\left( X,\ast \right) \diamond id_{Bin\left(
X\right) }$ and $\left( X,\circ \right) =O\left( X,\circ \right) \diamond
J\left( X,\circ \right) =\left( X,\ast \right) \diamond \left( X,\bullet
\right) $. This clearly shows the results of Theorem 4.3.1.\medskip

\noindent \textbf{Theorem 4.3.3 }\textsl{The right-zero-semigroup on X is }$%
j $\textsl{-composite. \smallskip }

\noindent \textit{Proof}. Let $\left( X,\bullet \right) $ be the
right-zero-semigroup on $X$. Suppose that $\left( X,\ast \right) =O\left(
X,\bullet \right) $ and $\left( X,\circ \right) =J\left( X,\bullet \right) $%
. By applying Proposition 4.2.5, $\left( X,\bullet \right) $ is $j$-normal.
Thus, $\left( X,\bullet \right) =\left( X,\ast \right) \diamond \left(
X,\circ \right) =\left( X,\circ \right) \diamond \left( X,\ast \right) $.
Consider $\left( X,\ast \right) $: (i)$\overline{\text{ }d^{\ast }}=%
\widetilde{D^{\diamond }}$, (ii) $x\ast y=x$, otherwise; and for $\left(
X,\circ \right) $: (i) $\overline{d^{\circ }}=\widetilde{d^{\bullet }}$,
(ii) $x\circ y=x\bullet y=y$, otherwise. Since neither one of the factors is
the left-zero-semigroup for $Bin\left( X\right) $, $\left( X,\bullet \right) 
$ is $j$-composite.

\begin{flushright}
$\blacksquare $
\end{flushright}

\noindent \textbf{Example 4.3.4} Let $\left( X,\bullet \right) $ be the
right-zero-semigroup as in Example 3.2.9 where $X=\{a,b,c\}$. Let $\left(
X,\ast \right) =O\left( X,\bullet \right) $ and $\left( X,\circ \right)
=J\left( X,\bullet \right) $, we can check that $\left( X,\bullet \right) $
is in fact $OJ$\textit{-} and $JO$-composite. Hence, $\left( X,\bullet
\right) $ is $j$-composite: 
\begin{equation*}
\text{\quad }%
\begin{tabular}{l|lll}
$\ast $ & $a$ & $b$ & $c$ \\ \hline
$a$ & $a$ & $a$ & $c$ \\ 
$b$ & $b$ & $b$ & $b$ \\ 
$c$ & $a$ & $c$ & $c$%
\end{tabular}%
\text{\quad }\diamond \text{\quad }%
\begin{tabular}{l|lll}
$\circ $ & $a$ & $b$ & $c$ \\ \hline
$a$ & $a$ & $b$ & $a$ \\ 
$b$ & $a$ & $b$ & $c$ \\ 
$c$ & $c$ & $b$ & $c$%
\end{tabular}%
\text{\quad }=\text{\quad }%
\begin{tabular}{l|lll}
$\bullet $ & $a$ & $b$ & $c$ \\ \hline
$a$ & $a$ & $b$ & $c$ \\ 
$b$ & $a$ & $b$ & $c$ \\ 
$c$ & $a$ & $b$ & $c$%
\end{tabular}%
\text{\quad }
\end{equation*}%
\noindent Moreover, its \textit{orient}-factor $\left( X,\ast \right) $ has
the following subtables:%
\begin{equation*}
\text{\quad 
\begin{tabular}{r|rr}
$\ast $ & $a$ & $b$ \\ \hline
$a$ & $a$ & $a$ \\ 
$b$ & $b$ & $b$%
\end{tabular}%
\quad 
\begin{tabular}{r|rr}
$\ast $ & $a$ & $c$ \\ \hline
$a$ & $a$ & $c$ \\ 
$c$ & $a$ & $c$%
\end{tabular}%
\quad 
\begin{tabular}{r|rr}
$\ast $ & $b$ & $c$ \\ \hline
$b$ & $b$ & $b$ \\ 
$c$ & $c$ & $c$%
\end{tabular}%
\quad }
\end{equation*}%
\noindent which implies that $\left( X,\ast \right) $ is
locally-zero.\bigskip

Given two distinct groupoids $\left( X,\triangleright \right) $ and $\left(
X,\triangleleft \right) $ in $Bin\left( X\right) $. Suppose that $\left(
X,\triangleright \right) \neq id_{Bin\left( X\right) }$ and $\left(
X,\triangleleft \right) \neq $ $id_{Bin\left( X\right) }$. Let $\left(
X,\bullet \right) $ be a groupoid such that $\left( X,\triangleright \right)
\neq \left( X,\bullet \right) \neq \left( X,\triangleleft \right) $. Then $%
\left( X,\bullet \right) $ is said to be:

\qquad (i) \emph{partially-right-prime, }$\partial _{r}$-prime, if $\left(
X,\bullet \right) =\left( X,\bullet \right) \diamond \left( X,\triangleright
\right) $;

\qquad (ii) \emph{partially-left-prime, }$\partial _{l}$-prime, if $\left(
X,\bullet \right) =\left( X,\triangleleft \right) \diamond \left( X,\bullet
\right) $.

Whence $\left( X,\triangleright \right) $ and $\left( X,\triangleleft
\right) $ behave like \textit{right}- and \textit{left}-identities
respectively. Here, $\left( X,\triangleright \right) $ and $\left(
X,\triangleleft \right) $ could be either $O\left( X,\bullet \right) $, $%
J\left( X,\bullet \right) $, $U\left( X,\bullet \right) $, $A\left(
X,\bullet \right) $ or any other factor of $\left( X,\bullet \right) $. The
next proposition demonstrates one such case.\bigskip

\noindent \textbf{Proposition 4.3.5} \textsl{A bi-diagonal groupoid is
partially-left-prime.\smallskip }

\noindent \textit{Proof.} Given a bi-diagonal groupoid $\left( X,\bullet
\right) $, then its \textit{skew}-factor $J(X,\bullet )=\left( X,\bullet
\right) $ since $\overline{d^{\circ }}=\widetilde{d^{\bullet }}=\overline{%
d^{\bullet }}$ and $x\circ y=x\bullet y$ otherwise. Meanwhile, its \textit{%
orient}-factor $O(X,\bullet )$ is not affected by the bi-diagonal property.
By Theorem 4.1.2, $(X,\bullet )$ has an $OJ$-factorization, 
\begin{eqnarray*}
(X,\bullet ) &=&O(X,\bullet )\diamond J(X,\bullet ) \\
&=&O(X,\bullet )\diamond (X,\bullet )\text{.}
\end{eqnarray*}

\noindent Therefore, $O(X,\bullet )$ is a left-identity in $(Bin\left(
X\right) ,\diamond )$ and the result follows.

\begin{flushright}
$\blacksquare $
\end{flushright}

\noindent \textbf{Example 4.3.6 }Consider the group $\left( X,\bullet
,e\right) $\ as defined in Example 4.5. Then clearly $\left( X,\bullet
,e\right) $ is bi-diagonal. Recall its \textit{orient-}factor $\left( X,\ast
,e\right) =O(X,\bullet ,e)$ and derive its \textit{skew}-factor $\left(
X,\circ ,e\right) =J(X,\bullet ,e)$ to obtain:%
\begin{equation*}
\quad 
\begin{tabular}{c|cccc}
$\ast $ & $e$ & $a$ & $b$ & $c$ \\ \hline
$e$ & $e$ & $e$ & $e$ & $c$ \\ 
$a$ & $a$ & $a$ & $b$ & $a$ \\ 
$b$ & $b$ & $a$ & $b$ & $b$ \\ 
$c$ & $e$ & $c$ & $c$ & $c$%
\end{tabular}%
\quad \text{and}\quad 
\begin{tabular}{c|cccc}
$\circ $ & $e$ & $a$ & $b$ & $c$ \\ \hline
$e$ & $e$ & $a$ & $b$ & $c$ \\ 
$a$ & $a$ & $e$ & $c$ & $b$ \\ 
$b$ & $b$ & $c$ & $a$ & $e$ \\ 
$c$ & $c$ & $b$ & $e$ & $a$%
\end{tabular}%
\quad
\end{equation*}

\noindent Then $(X,\bullet ,e)=O(X,\bullet ,e)\diamond (X,\bullet ,e)$ and
therefore the group $(X,\bullet ,e)$ is $\partial _{l}$-prime.

\section{\protect\bigskip Application}

Recall some of the algebras described in \textquotedblleft \textbf{Figure 1}%
\textquotedblright\ of Section 2.\medskip

We shall say an algebra $(X,\bullet ,0)$ of type $\left( 2,0\right) $ is a 
\emph{strong }$B1$\emph{-algebra} if it satisfies (B1) and equation \ref%
{strong}. Meaning, if for all $x,y\in X$,

\qquad (i) $x\bullet x=0$,

\qquad (ii) $x\bullet y=y\bullet x$ implies $x=y$.\medskip

A groupoid $\left( X,\bullet ,0\right) $ is \emph{semi-neutral} if for all $%
x,y\in X$,

\qquad (i) $x\bullet x=0$,

\qquad (ii) $x\bullet y=x$.\medskip

A $B1$-algebra $\left( X,\bullet ,0\right) $ is \emph{semi-neutral} if for $%
x\neq y$, $x\bullet y=x$ for all $x,y\in X.$\medskip

A normal/composite groupoid is \emph{semi-normal (resp., semi-composite)} if
only one of its factors is semi-neutral.\bigskip

\noindent \textbf{Proposition 5.1 }\textsl{A semi-neutral groupoid\ is
signature-prime and }$OJ$\textsl{-composite.}\smallskip

\noindent \textit{Proof}. Let $\left( X,\bullet ,0\right) $ be the
semi-neutral groupoid on $X$. Then $x\bullet y=x$ for all $x,y\in X$ and $%
x\bullet x=0$. Let $\left( X,\ast ,0\right) =U\left( X,\bullet ,0\right) $
and $\left( X,\circ ,0\right) =A\left( X,\bullet ,0\right) $, its \textit{%
signature}- and \textit{similar}-factors, respectively. Deriving them
according to \ref{Sig-Sim-F's} gives:%
\begin{equation*}
x\ast y=%
\begin{cases}
x & \text{ if }x=y, \\ 
x\bullet y=x & \text{otherwise.}%
\end{cases}%
\quad \text{and}\quad x\circ y=%
\begin{cases}
x\bullet x=0 & \text{ if }x=y, \\ 
x & \text{otherwise.}%
\end{cases}%
\end{equation*}%
\noindent Hence for all $x,y$ $\in X$, $\left( X,\bullet ,0\right)
=id_{Bin\left( X\right) }\diamond \left( X,\bullet ,0\right) $.

\noindent By Theorem 4.1.2, $\left( X,\bullet ,0\right) $ has an $OJ$%
-factorization. Let $\left( X,\circledast ,0\right) =O\left( X,\bullet
,0\right) $ and $\left( X,\odot ,0\right) =J\left( X,\bullet ,0\right) $,
its \textit{orient}- and \textit{skew}-factors, respectively. Deriving them
according to \ref{Orient-Skew-F's} gives: for $\circledast $: (i)$\overline{%
\text{ }d^{\circledast }}=\widetilde{D^{\diamond }}$, (ii) $x\circledast y=x$%
, otherwise; and for $\odot $: (i) $\overline{d^{\odot }}=\widetilde{%
d^{\bullet }}\neq D^{\diamond }$, (ii) $x\odot y=x\bullet y$, otherwise.
Thus, $\left( X,\ast ,0\right) \neq id_{Bin\left( X\right) }\neq \left(
X,\circ ,0\right) $ and $\left( X,\ast ,0\right) \neq \left( X,\bullet
,0\right) \neq \left( X,\circ ,0\right) $.

\begin{flushright}
$\blacksquare $
\end{flushright}

\noindent \textbf{Corollary 5.2} \textsl{A semi-neutral groupoid is
semi-normal.\smallskip }

\noindent \textit{Proof.} This is a direct result of Proposition 5.1 and the
definition of a semi-normal groupoid.

\begin{flushright}
$\blacksquare $
\end{flushright}

\noindent \textbf{Proposition 5.3} \textsl{The product of semi-neutral
groupoids is semi-neutral.}\vspace{0.07in}

\noindent \textit{Proof. }Consider semi-netural groupoids $\left( X,\ast
,0\right) $ and $\left( X,\circ ,0\right) $. Let $\left( X,\ast ,0\right)
\diamond \left( X,\circ ,0\right) =$ $\left( X,\bullet ,0\right) $ such that 
$x\bullet y=\left( x\ast y\right) \circ \left( y\ast x\right) $. Then, $%
x\bullet x=\left( x\ast x\right) \circ \left( x\ast x\right) =0$. If $x\neq
y,$.$x\bullet y=x\circ y$. It follows that $\left( X,\bullet ,0\right)
=\left( X,\circ ,0\right) $ and therefore is semi-neutral.

\begin{flushright}
$\blacksquare $
\end{flushright}

\noindent \textbf{Proposition 5.4} \textsl{The similar-factor of a }$B1$%
\textsl{-algebra is semi-neutral.\smallskip }

\noindent \textit{Proof.} Let $\left( X,\bullet ,0\right) $ be a $B1$%
-algebra. Consider the $AU$-factorization $\left( X,\bullet ,0\right)
=A\left( X,\bullet ,0\right) \diamond U\left( X,\bullet ,0\right) $. Let $%
\left( X,\ast ,0\right) :=U\left( X,\bullet ,0\right) $ and $\left( X,\circ
,0\right) :=A\left( X,\bullet ,0\right) $, its \textit{signature}- and 
\textit{similar}-factors, respectively. Deriving them according to \ref%
{Sig-Sim-F's} gives:%
\begin{equation*}
x\ast y=%
\begin{cases}
x & \text{ if }x=y, \\ 
x\bullet y & \text{otherwise.}%
\end{cases}%
\quad \text{and}\quad x\circ y=%
\begin{cases}
x\bullet x=0 & \text{ if }x=y, \\ 
x & \text{otherwise.}%
\end{cases}%
\end{equation*}%
\noindent Clearly, $\left( X,\circ ,0\right) $ is semi-neutral.

\begin{flushright}
$\blacksquare $
\end{flushright}

\noindent \textbf{Corollary 5.5} \textsl{A strong }$B1$\textsl{-algebra is
semi-normal.\smallskip }

\noindent \textit{Proof.} This is a direct result of Corollary 3.2.5,
Proposition 5.4 and the definition of a semi-normal algebra.

\begin{flushright}
$\blacksquare $
\end{flushright}

\noindent \textbf{Corollary 5.6} \textsl{A strong }$B1$\textsl{-algebra }$%
\left( X,\bullet ,0\right) $ \textsl{is semi-composite if it is not
semi-neutral, i.e., if }$x\bullet y\neq x$ for all $x,$ $y\in X$.\textsl{%
\smallskip }

\noindent \textit{Proof.} Let $\left( X,\bullet ,0\right) $ be a strong $B1$%
-algebra. Let $\left( X,\ast ,0\right) :=U\left( X,\bullet ,0\right) $ and $%
\left( X,\circ ,0\right) :=A\left( X,\bullet ,0\right) $, its \textit{%
signature}- and \textit{similar}-factors respectively. Deriving them
according to \ref{Sig-Sim-F's}. Assume that $x\bullet y=x$. Then $x\ast y=x$
for all $x,y\in X$. Thus, $\left( X,\bullet ,0\right) =id_{Bin\left(
X\right) }\diamond (X,\bullet ,0)$ which makes it signature-prime and not $u$%
-composite.

\begin{flushright}
$\blacksquare $
\end{flushright}

\noindent \textbf{Example 5.7} Let $\left( X,\bullet ,0\right) =\left(
\left\{ 0,1,2\right\} ,\bullet \right) $ be a strong $BCK$-algebra of order
3 where \textquotedblleft $\bullet $\textquotedblright\ is defined by the
following Cayley table:

\begin{center}
$%
\begin{tabular}{l|lll}
$\bullet $ & 0 & 1 & 2 \\ \hline
0 & 0 & 0 & 0 \\ 
1 & 1 & 0 & 1 \\ 
2 & 2 & 2 & 0%
\end{tabular}%
\ $
\end{center}

\noindent Let $\left( \left\{ 0,1,2\right\} ,\ast \right) =U\left( \left\{
0,1,2\right\} ,\bullet \right) $ and $\left( \left\{ 0,1,2\right\} ,\circ
\right) =A\left( \left\{ 0,1,2\right\} ,\bullet \right) $. Its $UA$%
.-factorization is:%
\begin{equation*}
\quad 
\begin{tabular}{l|lll}
$\ast $ & 0 & 1 & 2 \\ \hline
0 & 0 & 0 & 0 \\ 
1 & 1 & 1 & 1 \\ 
2 & 2 & 2 & 2%
\end{tabular}%
\quad \diamond \quad 
\begin{tabular}{l|lll}
$\circ $ & 0 & 1 & 2 \\ \hline
0 & 0 & 0 & 0 \\ 
1 & 1 & 0 & 1 \\ 
2 & 2 & 2 & 0%
\end{tabular}%
\quad =\quad 
\begin{tabular}{l|lll}
$\bullet $ & 0 & 1 & 2 \\ \hline
0 & 0 & 0 & 0 \\ 
1 & 1 & 0 & 1 \\ 
2 & 2 & 2 & 0%
\end{tabular}%
\quad 
\end{equation*}%
\noindent Therefore, $\left( \left\{ 0,1,2\right\} ,\bullet \right) $ is 
\textit{signature}-prime and $u$-normal. Moreover, $\left( \left\{
0,1,2\right\} ,\bullet \right) $ as defined is semi-neutral. Next, \textit{%
derive} its \textit{orient}- and \textit{skew}-factors $O\left( X,\bullet
,0\right) $ and $J\left( X,\bullet ,o\right) $, respectively. Let $\left(
X,\circledast ,0\right) =O\left( X,\bullet ,0\right) $ and $\left( X,\odot
,0\right) =J\left( X,\bullet ,0\right) $. We have the following product:%
\textsl{\smallskip }%
\begin{equation*}
\text{\quad }%
\begin{tabular}{l|lll}
$\ast $ & 0 & 1 & 2 \\ \hline
0 & 0 & 0 & 2 \\ 
1 & 1 & 1 & 1 \\ 
2 & 0 & 2 & 2%
\end{tabular}%
\text{\quad }\diamond \text{\quad }%
\begin{tabular}{l|lll}
$\circ $ & 0 & 1 & 2 \\ \hline
0 & 0 & 0 & 2 \\ 
1 & 1 & 0 & 1 \\ 
2 & 0 & 2 & 0%
\end{tabular}%
\text{\quad }=\text{\quad }%
\begin{tabular}{l|lll}
$\bullet $ & 0 & 1 & 2 \\ \hline
0 & 0 & 0 & 0 \\ 
1 & 1 & 0 & 1 \\ 
2 & 2 & 2 & 0%
\end{tabular}%
\end{equation*}%
\noindent Hence,  $O\left( X,\bullet ,0\right) \neq id_{Bin\left( X\right)
}\neq J\left( X,\bullet ,0\right) $ implies that $\left( X,\bullet ,0\right) 
$ is $OJ$-composite.\medskip 

\noindent \textbf{Example 5.8} Let $\left( X,\bullet ,0\right) =\left\{
(0,1,2),\bullet \right\} $ be a strong $Q$-algebra of order 3 where
\textquotedblleft $\bullet $\textquotedblright\ is given by the following
Cayley table:

\begin{center}
$%
\begin{tabular}{l|lll}
$\bullet $ & 0 & 1 & 2 \\ \hline
0 & 0 & 2 & 1 \\ 
1 & 1 & 0 & 2 \\ 
2 & 2 & 1 & 0%
\end{tabular}%
$
\end{center}

\noindent Let $\left( \left\{ 0,1,2\right\} ,\ast \right) =U\left( \left\{
0,1,2\right\} ,\bullet \right) $ and $\left( \left\{ 0,1,2\right\} ,\circ
\right) =A\left( \left\{ 0,1,2\right\} ,\bullet \right) $. Its $UA$%
.-factorization is:%
\begin{equation*}
\quad 
\begin{tabular}{l|lll}
$\ast $ & 0 & 1 & 2 \\ \hline
0 & 0 & 2 & 1 \\ 
1 & 1 & 1 & 2 \\ 
2 & 2 & 1 & 2%
\end{tabular}%
\quad \diamond \quad 
\begin{tabular}{l|lll}
$\circ $ & 0 & 1 & 2 \\ \hline
0 & 0 & 0 & 0 \\ 
1 & 1 & 0 & 1 \\ 
2 & 2 & 2 & 0%
\end{tabular}%
\quad =\quad 
\begin{tabular}{l|lll}
$\bullet $ & 0 & 1 & 2 \\ \hline
0 & 0 & 2 & 1 \\ 
1 & 1 & 0 & 2 \\ 
2 & 2 & 1 & 0%
\end{tabular}%
\quad 
\end{equation*}%
\noindent Since $\left( \left\{ 0,1,2\right\} ,\ast ,0\right) $ $\neq
Id_{Bin\left( X\right) }$ and $\left( \left\{ 0,1,2\right\} ,\circ ,0\right) 
$ is semi-neutral, we can conclude that $\left( \left\{ 0,1,2\right\}
,\bullet ,0\right) $ is semi-composite.\bigskip 

P.J. Allen, H.S. Kim and Neggers in \cite{AKN05} introduced the notion of
Smarandache disjointness in algebras. Two groupoids (algebras) $(X,\bullet )$
and $(X,\ast )$ are said to be Smarandache disjoint if we add some axioms of
an algebra $(X,\bullet )$ to an algebra $\left( X,\ast \right) $, then the
algebra $\left( X,\ast \right) $ becomes a trivial algebra, i.e., $%
\left\vert X\right\vert =1$.\bigskip 

\noindent \textbf{Proposition 5.9 }\textsl{The class of abelian groupoids
and the class of }$u$\textsl{-normal groupoids are Smarandache
disjoint.\smallskip }

\noindent \textit{Proof. } Let $\left( X,\bullet \right) $ $\in Ab(X)$, the
collection of all abelian groupoids defined on $X$. Suppose that $(X,\circ
)=A(X,\bullet )$ and $(X,\ast )=U$ $\left( X,\bullet \right) $. By Theorem
3.2.3, $\left( X\bullet \right) $ admits an $AU$-factorization. Consider $%
(X,\ast )\diamond (X,\circ )$, then for $x=y$, 
\begin{eqnarray*}
x\diamond x &=&\left( x\ast x\right) \circ \left( x\ast x\right)  \\
&=&x\circ x \\
&=&x\bullet x.
\end{eqnarray*}%
If $x\neq y$, 
\begin{eqnarray*}
x\diamond y &=&\left( x\ast y\right) \circ \left( y\ast x\right)  \\
&=&\left( x\bullet y\right) \circ \left( y\bullet x\right)  \\
&=&\left( x\bullet y\right) \bullet \left( x\bullet y\right) .
\end{eqnarray*}%
\noindent Hence, $\left( X\bullet \right) $ admits a $UA$-factorization only
if $\left( x\bullet y\right) \bullet \left( x\bullet y\right) =x\bullet y.$
This means that $\left( X,\bullet \right) $ is $u$-normal if it is either
the left- or right-zero-semigroup. Since both such groupoids are not
abelian, then $X$ must only have one element and the conclusion follows.

\begin{flushright}
$\blacksquare $
\end{flushright}

Suppose that in $Bin(X)$ we consider all those groupoids $(X,\ast )$ with
the orientation property. Thus, $x\ast x=x$ as a consequence. If $(X,\ast )$
and $(X,\circ )$ both have the orientation property, then for $x\diamond
\,\,y=(x\ast y)\circ (y\ast x)$ we have the possibilities: $x\ast
x=x,\,y\ast y=y,\,x\ast y\in \{x,y\}$ and $y\ast x\in \{x,y\}$, so that $%
x\,\diamond y\in \{x,y\}$. It follows that if $OP(X)$ denotes this
collection of groupoids, then $(OP(X),\diamond )$ is a subsemigroup \cite%
{KN08} of $(Bin(X),$ $\diamond )$.

In a sequence of papers Nebesk\'{y} (\cite{Ne98}, \cite{Ne00}, \cite{Ne06})
associated with graphs $(V,E)$ groupoids $(V,\ast )$ with various properties
and conversely. He defined a \textit{travel groupoid} $(X,\ast )$ as a
groupoid satisfying the axioms: $(u\ast v)\ast u=u$ and $(u\ast v)\ast v=u$
implies $u=v$. If one adds these two laws to the orientation property, then $%
(X,\ast )$ is an OP-travel-groupoid. In this case $u\ast v=v$ implies $v\ast
u=u$, i.e., $uv\in E$ implies $vu\in E$, i.e., the digraph $(X,E)$ is a
(simple) graph if $uu\not\in E$, with $u\ast u=u$. Also, if $u\not=v$, then $%
u\ast v=u$ implies $(u\ast v)\ast v=u\ast v=u$ is impossible, whence $u\ast
v=v$ and $uv\in E$, so that $(X,E)$ is a complete (simple) graph. \medskip

In a recent paper, Ahn, Kim and Neggers \cite{AKN18} related graphs with
binary systems in the center of $Bin\left( X\right) $. Given an element of $%
ZBin(X)$, say $(X,\bullet )$ , they constructed a graph, $\Gamma _{X}$ by
letting $V(\Gamma _{X})=X$ and $(x,y)\in E(\Gamma _{X})$, the edge set of $%
\Gamma _{X}$, such that $x\neq y$, $y\bullet x=y$ and $x\bullet y=x$. Thus,
if $(x,y)\in E(\Gamma _{X})$, then $(y,x)\in E(\Gamma _{X})$ as well and
they identify $(x,y)=(y,x)$ as an undirected edge of $\Gamma _{X}$. Then
they concluded that if $\left( X,\bullet \right) $ is the
left-zero-semigroup, then $\Gamma _{X}$ is the complete graph on $X$. Also,
if $\left( X,\bullet \right) $ is the right-zero-semigroup, then $\Gamma
_{X} $ is the null graph on $X$, since $E(\Gamma _{X})=\varnothing $.
\bigskip

\noindent \textbf{Example 5.10} Let $X=\left\{ a,b,c,d\right\} $ and
consider the simple graph on $X$: 
\begin{equation*}
\text{\begin{tikzpicture} \node[vertex][fill] (a) at (1,2)
[label=above:$a$]{}; \node[vertex][fill] (b) at (0,1) [label=left:$b$] {};
\node[vertex][fill] (c) at (2,1) [label=right:$c$] {}; \node[vertex][fill]
(d) at (1,0) [label=below:$d$] {}; \path (a) edge (b) (d) edge (b) (c) edge
(b) ; \end{tikzpicture}}
\end{equation*}%
\noindent Then the associated groupoid table with binary operation
\textquotedblleft $\bullet $\textquotedblright\ is:

\begin{equation*}
\begin{tabular}{c|cccc}
$\bullet $ & a & b & c & d \\ \hline
a & a & a & c & d \\ 
b & b & b & b & b \\ 
c & a & c & c & d \\ 
d & a & d & c & d%
\end{tabular}%
\end{equation*}%
\bigskip

\noindent By applying Proposition 4.2.5 to $\left( X,\bullet \right) $, we
have the product of $O\left( X,\bullet \right) $ and $J\left( X,\bullet
\right) $ given by their respective tables:

\begin{equation*}
\text{\quad }%
\begin{tabular}{c|cccc}
$\ast $ & a & b & c & d \\ \hline
a & a & a & a & d \\ 
b & b & b & c & b \\ 
c & c & b & c & c \\ 
d & a & d & d & d%
\end{tabular}%
\text{\quad }\diamond \text{\quad }%
\begin{tabular}{c|cccc}
$\circ $ & a & b & c & d \\ \hline
a & a & a & c & a \\ 
b & b & b & c & b \\ 
c & a & b & c & d \\ 
d & d & d & c & d%
\end{tabular}%
\text{\quad }=\text{\quad }%
\begin{tabular}{c|cccc}
$\bullet $ & a & b & c & d \\ \hline
a & a & a & c & d \\ 
b & b & b & b & b \\ 
c & a & c & c & d \\ 
d & a & d & c & d%
\end{tabular}%
\text{\quad }
\end{equation*}%
\bigskip

\noindent We can visualize this product with the associated graphs of
groupoids $\left( X,\ast \right) $ and $\left( X,\circ \right) $:

\begin{equation*}
\text{\quad }%
\begin{tabular}{cc}
\text{\begin{tikzpicture} \node[vertex][fill] (a) at (1,2)[label=above:$a$]
{}; \node[vertex][fill] (b) at (0,1) [label=above:$b$] {};
\node[vertex][fill] (c) at (2,1) [label=above:$c$] {}; \node[vertex][fill]
(d) at (1,0) [label=below:$d$] {}; \path (a) edge (b) (a) edge (c) (b) edge
(d) (c) edge (d) ; \end{tikzpicture}} & 
\end{tabular}%
\diamond \text{\quad }%
\begin{tabular}{cc}
\text{\begin{tikzpicture} \node[vertex][fill] (a) at (1,2) [label=above:$a$]
{}; \node[vertex][fill] (b) at (0,1) [label=above:$b$] {};
\node[vertex][fill] (c) at (2,1) [label=above:$c$] {}; \node[vertex][fill]
(d) at (1,0) [label=below:$d$] {}; \path (a) edge (b) (a) edge (d) (b) edge
(d) ; \end{tikzpicture}} & 
\end{tabular}%
\text{\quad }=\text{\quad }%
\begin{tabular}{cc}
\text{\begin{tikzpicture} \node[vertex][fill] (a) at (1,2) [label=above:$a$]
{}; \node[vertex][fill] (b) at (0,1) [label=above:$b$] {};
\node[vertex][fill] (c) at (2,1) [label=above:$c$] {}; \node[vertex][fill]
(d) at (1,0) [label=below:$d$] {}; \path (a) edge (b) (d) edge (b) (c) edge
(b) ; \end{tikzpicture}} & 
\end{tabular}%
\text{\quad }
\end{equation*}

Thus, any simple graph constructed in this manner can be decomposed into two
or more other factors with the binary product \textquotedblleft $\diamond $%
\textquotedblright . This fact is further illustrated in the next
example.\medskip

\noindent \textbf{Example 5.11 }Let $\left( X,\bullet \right) =\left(
\left\{ 0,1,2,3,4,5\right\} ,\bullet \right) $ be the locally-zero groupoid
defined as in Example 4.3.2. Then its associated graph decomposes into its
factors $\left( X,\ast \right) $ and $\left( X,\circ \right) $:

\begin{equation*}
\text{\quad }%
\begin{tabular}{cc}
\text{\begin{tikzpicture} \node[vertex][fill] (0) at (1,2) [label=above:$0$]
{}; \node[vertex][fill] (1) at (0,1) [label=above:$1$] {};
\node[vertex][fill] (2) at (1,0) [label=below:$2$] {}; \node[vertex][fill]
(3) at (2,0) [label=below:$3$] {}; \node[vertex][fill] (4) at (3,1)
[label=above:$4$] {}; \node[vertex][fill] (5) at (2,2) [label=above:$5$] {};
\path (0) edge (2) (0) edge (3) (0) edge (5) (1) edge (4) (2) edge (5) (3)
edge (4) (3) edge (5) (4) edge (5) ; \end{tikzpicture}} & 
\end{tabular}%
\text{\quad }=\text{\quad }%
\begin{tabular}{cc}
\text{\begin{tikzpicture} \node[vertex][fill] (0) at (1,2) [label=above:$0$]
{}; \node[vertex][fill] (1) at (0,1) [label=above:$1$] {};
\node[vertex][fill] (2) at (1,0) [label=below:$2$] {}; \node[vertex][fill]
(3) at (2,0) [label=below:$3$] {}; \node[vertex][fill] (4) at (3,1)
[label=above:$4$] {}; \node[vertex][fill] (5) at (2,2) [label=above:$5$] {};
\path (0) edge (2) (0) edge (1) (0) edge (3) (0) edge (4) (1) edge (2) (1)
edge (3) (1) edge (5) (2) edge (4) (2) edge (5) (3) edge (4) (3) edge (5)
(4) edge (5) ; \end{tikzpicture}} & 
\end{tabular}%
\text{\quad }\diamond \text{\quad }%
\begin{tabular}{cc}
\text{\begin{tikzpicture} \node[vertex][fill] (0) at (1,2) [label=above:$0$]
{}; \node[vertex][fill] (1) at (0,1) [label=above:$1$] {};
\node[vertex][fill] (2) at (1,0) [label=below:$2$] {}; \node[vertex][fill]
(3) at (2,0) [label=below:$3$] {}; \node[vertex][fill] (4) at (3,1)
[label=above:$4$] {}; \node[vertex][fill] (5) at (2,2) [label=above:$5$] {};
\path (0) edge (2) (0) edge (3) (2) edge (3) (2) edge (5) (3) edge (4) (3)
edge (5) (4) edge (5) ; \end{tikzpicture}} & 
\end{tabular}%
\text{\quad }
\end{equation*}

\section{\protect\bigskip Generalization and Summary}

In this final note, we discuss two generalizations which can serve as
grounds for future exploration of groupoid factorizations or algebra
decompositions via the groupoid product \textquotedblleft $\diamond $%
\textquotedblright .

\subsection{$\Psi $\textit{-type-Factorization}}

Let $\Psi $ be a groupoid operation that interchanges elements of any two
given groupoids and produces two other (possibly identical) groupoids. Given
groupoid $\left( X,\bullet \right) \in Bin\left( X\right) $ and the
left-zero-semigroup as $id_{Bin\left( X\right) }$, define $\Psi :Bin\left(
X\right) \times Bin\left( X\right) \rightarrow Bin\left( X\right) \times
Bin\left( X\right) $. A $\Psi $\textit{-\emph{type-factorization} }of $%
\left( X,\bullet \right) $ gives a pair of groupoid factors as follows:%
\begin{equation*}
\Psi (\left( X,\bullet \right) ,id_{Bin\left( X\right) })=\left( \left(
X,\bullet \right) _{L},\left( X,\bullet \right) _{R}\right)
\end{equation*}%
where $\left( X,\bullet \right) _{L}=\Psi _{\alpha }(\left( X,\bullet
\right) ,id_{Bin\left( X\right) })$ and $\left( X,\bullet \right) _{R}=\Psi
_{\alpha }(id_{Bin\left( X\right) },\left( X,\bullet \right) )$, the\emph{\
left- }and\emph{\ right-}$\Psi $\emph{-factors} of $\left( X,\bullet \right) 
$, respectively, such that the maps $\Psi _{\alpha }$ and $\alpha $ are
defined as $\Psi _{\alpha }:Bin\left( X\right) \times Bin\left( X\right)
\rightarrow Bin\left( X\right) $ and $\alpha :Bin\left( X\right) \rightarrow
Bin\left( X\right) $.

Let $\left( X,\ast \right) :=\left( X,\bullet \right) _{L}$ and $\left(
X,\circ \right) :=\left( X,\bullet \right) _{R}$, then $\left( X,\bullet
\right) $ can be represented as a product of the groupoid pair, i.e., 
\begin{eqnarray*}
\left( X,\bullet \right) &=&\left( X,\ast \right) \diamond \left( X,\circ
\right) \text{ and/or } \\
\left( X,\bullet \right) &=&\left( X,\circ \right) \diamond \left( X,\ast
\right)
\end{eqnarray*}

thus rendering $\left( X,\bullet \right) $ as:

\qquad (i) $\Psi $-\emph{prime}, if $\left( X,\bullet \right)
_{L}=id_{Bin\left( X\right) }$ or $\left( X,\bullet \right)
_{R}=id_{Bin\left( X\right) }$; or

\qquad (ii) $\Psi $-\emph{normal} if $\left( X,\ast \right) \diamond \left(
X,\circ \right) =\left( X,\circ \right) \diamond \left( X,\ast \right) $; or

\qquad (iii) $\Psi $-\emph{composite} if $\left( X,\bullet \right) $ is $%
\Psi $-normal but not $\Psi $-prime.\medskip

An example of this $\Psi $-type-factorization is our first method of \textit{%
similar}-\textit{signature}-factorization where 
\begin{equation*}
\Psi _{d}(\left( X,\bullet \right) ,id_{Bin\left( X\right) })=\{\left(
X,\bullet \right) |d\left( X,\bullet \right) =d\left( id_{Bin\left( X\right)
}\right) \}
\end{equation*}%
and 
\begin{equation*}
\Psi _{d}(id_{Bin\left( X\right) },\left( X,\bullet \right)
)=\{id_{Bin\left( X\right) }|d\left( id_{Bin\left( X\right) }\right)
=d\left( X,\bullet \right) \}
\end{equation*}%
The $\Psi $ in that case switched the diagonal $d$ of the parent groupoid $%
\left( X,\bullet \right) $ with that of the left-zero-semigroup, $%
id_{Bin\left( X\right) }$, to obtain the \textit{signature}- and \textit{%
similar}-factors $\left( X,\circ \right) $ and $\left( X,\ast \right) $,
respectively. Hence, the \textit{signature-} and \textit{similar}-factors of
a groupoid are $\Psi $\textit{-}\emph{type-factors}\textit{.}

\subsection{$\protect\tau $\textit{-type-Factorization}}

Let $\tau $ be a groupoid operation that manipulates elements of any given
pair of groupoid in the same fashion. Given groupoid $\left( X,\bullet
\right) \in Bin\left( X\right) $ and the left-zero-semigroup as $%
id_{Bin\left( X\right) }$, define $\tau :Bin\left( X\right) \times Bin\left(
X\right) \rightarrow Bin\left( X\right) \times Bin\left( X\right) $. A $\tau 
$\textit{-}\emph{type-factorization} of $\left( X,\bullet \right) $ is given
as follows:%
\begin{equation*}
\tau \left( \left( X,\bullet \right) ,id_{Bin\left( X\right) }\right)
=\left( \left( X,\bullet \right) _{L},\left( X,\bullet \right) _{R}\right)
\end{equation*}%
where $\left( X,\bullet \right) _{L}=\theta (id_{Bin\left( X\right) })$ and
and $\left( X,\bullet \right) _{R}=\theta \left( X,\bullet \right) $ such
that the map $\theta :Bin\left( X\right) \rightarrow Bin\left( X\right) $,
the \emph{left- }and\emph{\ right-}$\tau $\emph{-factors}\textit{\ }of $%
\left( X,\bullet \right) $, respectively. Let $\left( X,\ast \right)
:=\left( X,\bullet \right) _{L}$ and $\left( X,\circ \right) :=\left(
X,\bullet \right) _{R}$, then $\left( X,\bullet \right) $ could factor into
a product of the groupoid pair, i.e., 
\begin{eqnarray*}
\left( X,\bullet \right) &=&\left( X,\ast \right) \diamond \left( X,\circ
\right) \text{ and/or } \\
\left( X,\bullet \right) &=&\left( X,\circ \right) \diamond \left( X,\ast
\right) \text{.}
\end{eqnarray*}%
Once again rendering $\left( X,\bullet \right) $ as:

\qquad (i) $\tau $-\emph{prime}, if $\left( X,\bullet \right)
_{L}=id_{Bin\left( X\right) }$ or $\left( X,\bullet \right)
_{R}=id_{Bin\left( X\right) }$; or

\qquad (ii) $\tau $-\emph{normal} if $\left( X,\ast \right) \diamond \left(
X,\circ \right) =\left( X,\circ \right) \diamond \left( X,\ast \right) $; or

\qquad (iii) $\tau $-\emph{composite} if $\left( X,\bullet \right) $ is $%
\tau $-normal but not $\tau $-prime.\medskip

An example of this $\tau $-type-factorization is our second method of 
\textit{orient-skew-}factorization where $O\left( X,\bullet \right) :=\left(
X,\bullet \right) _{L}$ and $J\left( X,\bullet \right) :=\left( X,\bullet
\right) _{R}$. The $\tau $ (indeed, $\theta $) in that scenario reversed the
anti-diagonal of a given groupoid. Hence, applying $\tau $ to the
left-zero-semigroup $id_{Bin\left( X\right) }$ and to the parent groupoid $%
\left( X,\bullet \right) $ results in the \textit{orient}- and \textit{skew}%
-factors $\left( X,\ast \right) $ and $\left( X,\circ \right) $,
respectively. In conclusion, the \textit{orient-} and \textit{skew}-factors
of a groupoid are $\tau $\emph{-type-factors}\textit{.}

\subsection{Summary}

The goal of this paper was to gain more insight about the dynamics of binary
systems, namely groupoids or algebras equipped with a single binary
operation. We have shown that a strong groupoid can be represented as a
\textquotedblleft composite\textquotedblright\ groupoid of its \textit{%
similar-} and \textit{signature-} derived factors. Moreover, we concluded
that an idempotent groupoid with the orientation property, can be decomposed
into a product of its \textit{orient-} and \textit{skew-} factors. An
application into the fields of logic-algebras and graph theory were briefly
introduced. We found that a \textit{semi}-neutral $B1$-algebra is \textit{%
signature}-prime, $OJ$-composite and \textit{semi}-normal. Meanwhile, a
strong $B1$-algebra is then \textit{semi}-composite if it is not \textit{semi%
}-neutral. We finished our note with generalizations of our two methods in
hopes that other factorizations can be discovered in the near future. It may
be interesting to find other conditions for a groupoid to have such
decompositions. As a final reminder, factorization can be useful in various
applications such as algebraic cryptography and DNA code theory. We intend
to extend our investigation in the future to hypergroupoid, semigroups as
well as determine other factorizations and explore their applications.

\bigskip

\section{Acknowledgment}

The author is grateful to J. Neggers, H.S. Kim, C. Odenthal and the referee
for their valuable suggestions and help.

This research did not receive any specific grant from funding agencies in
the public, commercial, or not-for-profit sectors.

\end{document}